%% file: main.tex
\documentclass[a4paper]{svmult}

\input{custom_commands}

\title*{Properties of Spline Spaces Over Structured Hierarchical Box Partitions}
\author{Ivar Stangeby \and Tor Dokken}
\institute{
  Ivar Stangeby \at SINTEF, Forskningsveien 1, 0373 Oslo, Norway
  \protect \email{ivar.stangeby@sintef.no}
  \and
  Tor Dokken \at SINTEF, Forskningsveien 1, 0373 Oslo, Norway
  \protect \email{tor.dokken@sintef.no}
}

\begin{document}

\maketitle

\abstract{
  Given a spline space spanned by Truncated Hierarchical B-splines (THB), it is
  always possible to construct a spline space spanned by Locally Refined
  B-splines (LRB) that contains the THB-space. Starting from configurations
  where the two spline spaces are equal, we adress what happens to the
  properties of the LRB-space when it is modified by local one-directional
  refinement at convex corners of, and along edges between dyadic refinement
  regions. We show that such local modifications can reduce the number of
  B-splines over each element to the minimum prescribed by the polynomial
  bi-degree, and that such local refinements can be used for improving the
  condition numbers of mass and stiffness matrices.
}

\input{introduction}

\input{boundary_multiplicity}
\input{box_partitions_meshes_spline_spaces}

\input{local_modifications_and_overloading}
\input{numerical_experiments}

\input{conclusion}

\section{Acknowledgement}
This project has received funding from the The Research Council of Norway under grant agreement No 270922.

\clearpage
\bibliography{bibliography} 
\bibliographystyle{ieeetr}

\end{document}

%% file: custom_commands.tex
\usepackage{microtype}

\usepackage{newtxtext}
\usepackage{newtxmath}
\usepackage{makeidx} 

\usepackage{enumitem}
\setlist[description]{before = \leavevmode, font=\normalfont\bfseries}
\setlist[itemize]{font = \upshape, before = \leavevmode}
\setlist[enumerate]{font = \upshape, before = \leavevmode}

\usepackage[section]{algorithm}
\usepackage{algpseudocode}

\usepackage[]{booktabs} 

\usepackage[]{helvet} 
\usepackage[]{courier} 

\usepackage[]{commath}
\usepackage[]{mathtools}

\usepackage{graphicx}
\graphicspath{{images/}}
\usepackage{subfig}

\usepackage{varioref}
\usepackage[nameinlink, capitalize, noabbrev]{cleveref}

\newcommand{\define}[1]{\textbf{#1}}

\newcommand{\boxpartition}{\ensuremath{\mathcal{E}}}
\newcommand{\dbox}{\ensuremath{\beta}}
\newcommand{\element}{\ensuremath{\beta}}
\newcommand{\mesh}{\ensuremath{\mathcal{M}}}
\newcommand{\meshline}{\ensuremath{\gamma}}

\newcommand{\bspline}{\ensuremath{B[\knotvector^1, \ldots, \knotvector^d]}}
\newcommand{\bsplines}{\ensuremath{\mathcal{B}}}

\newcommand{\thbsplines}{\ensuremath{\mathcal{H}}}
\newcommand{\lrsplines}{\ensuremath{\mathcal{L}}}

\newcommand{\splinespace}[3]{\ensuremath{\mathcal{S}_{#2}^{#3}(#1)}}
\newcommand{\continuous}[1]{\ensuremath{C^{#1}}}
\newcommand{\polynomials}{\ensuremath{\Pi}}

\newcommand{\R}{\ensuremath{\mathbb{R}}}
\newcommand{\N}{\ensuremath{\mathbb{N}}}
\newcommand{\interior}[1]{%
  {\kern0pt#1}^{\mathrm{o}}%
}
\newcommand{\closure}[1]{%
  \overline{#1}
}

\newcommand{\knot}{\ensuremath{t}}
\newcommand{\knotvector}{\ensuremath{\vec{\knot}}}
\newcommand{\Supp}{\ensuremath{\mathrm{supp}}}
\newcommand{\boundary}{\ensuremath{\partial}}
\newcommand{\SPAN}{\ensuremath{\mathrm{span}}}
\newcommand{\trunc}{\ensuremath{\mathrm{trunc}}}

\newcommand{\mat}[1]{\ensuremath{\vec{#1}}}
\newcommand{\laplace}[1]{\ensuremath{\Delta #1}}
\newcommand{\grad}[1]{\ensuremath{\nabla#1}}
\newcommand{\cond}[1]{\ensuremath{\mathrm{Cond(#1)}}}

%% file: introduction.tex
\section{Introduction}

The use of Hierarchical B-splines (HB) introduced in
\cite{forseyHierarchicalBsplineRefinement1988} has gained much attention in
Isogeometric Analysis (IgA) in recent years. Hierarchical B-splines are based
on a dyadic sequence of grids determined by scaled lattices. On each
hierarchical level a spline space is defined as the tensor product of
univariate spline spaces spanned by uniform B-splines.

Hierarchical B-splines do not constitute a partition of unity, a much desired
property in both Computer Aided Design (CAD) and IgA. As a remedy to this 
Truncated Hierarchical B-splines (THB) \cite{THBsplinesTruncatedBasis2012,
VUONG20113554} were introduced, where B-splines on one hierarchical level are
suitably \emph{truncated} by B-splines from finer hierarchical levels when
the support of a B-spline at a finer level is contained in the support of a B-spline at a coarser level. 

An alternative to the THB-approach for forming a partition of unity came with
the introduction of Locally Refined B-splines (LRB)
\cite{dokkenPolynomialSplinesLocally2013}, where initial tensor product
B-splines are split until only B-splines of minimal support remain.  LRB
permits dyadic refinement of hierarchical meshes while ensuring that all
B-splines have minimal support. In the occasional case where meshlines at a
dyadic level are too short to split an LR B-spline, the meshlines in question
are extended. This fact ensures that the spline space spanned by THB-splines is
either identical to or constitutes a subset of the LRB spline space.

\begin{figure}[htbp]
  \centering
  \subfloat[]{\includegraphics[width=0.49\linewidth]{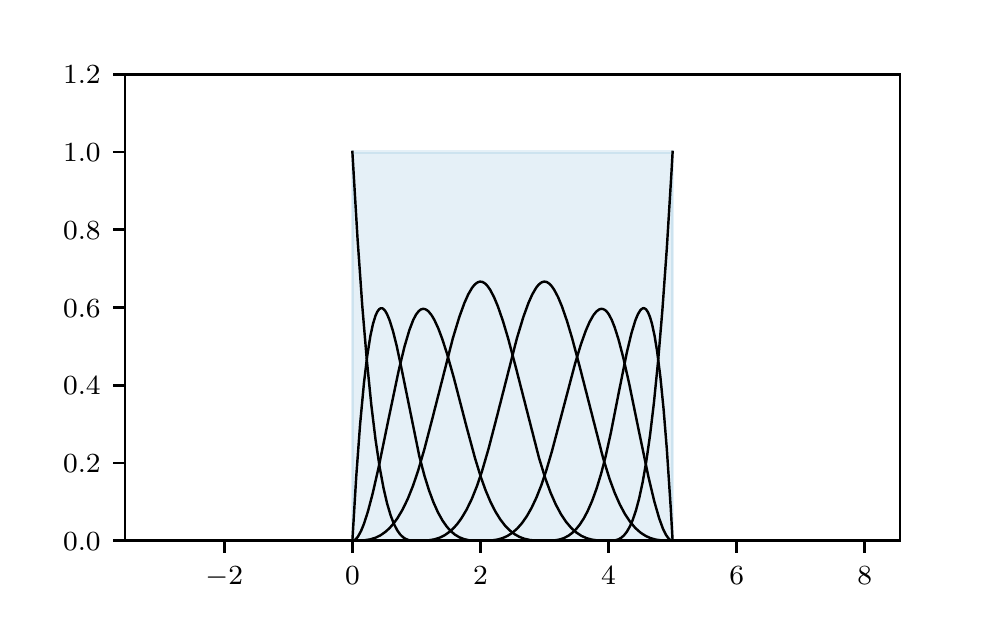}\label{sub:open}}
  \subfloat[]{\includegraphics[width=0.49\linewidth]{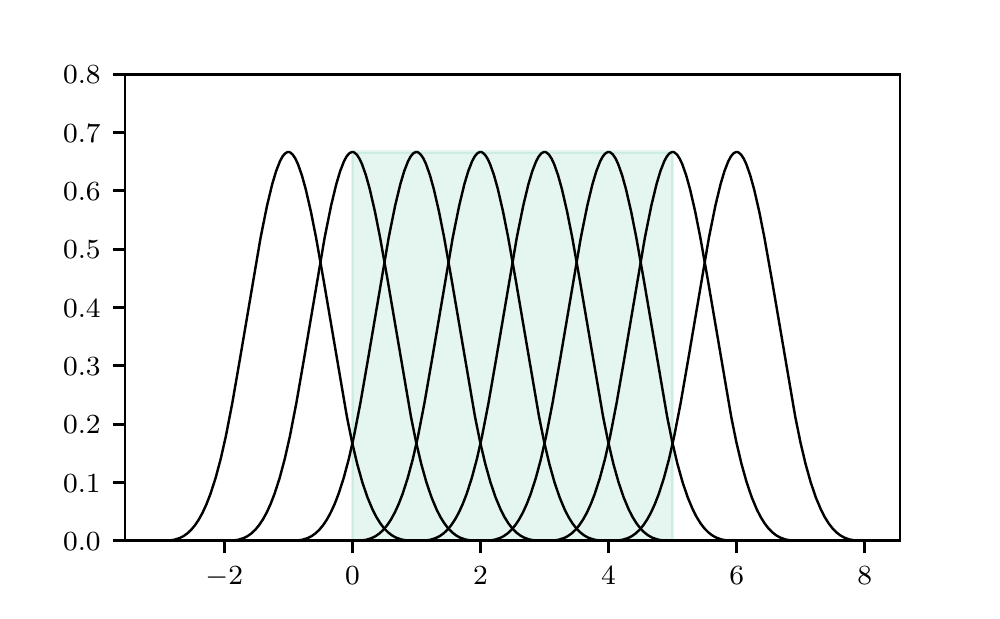}\label{sub:ghost}}
  \caption[Ghost domain and open knot vectors.]{Spline spaces over the domain \( \Omega = [0, 5] \). 
    In \protect\subref{sub:open}, the partition of unity is satisfied at the boundary
    by setting the knot multiplicity to \( m = d + 1 = 4 \). In
    \protect\subref{sub:ghost}, the partition of unity is satisfied at the boundary by
    extending the domain to allow the full polynomial space to be spanned at
    the boundary elements. The shaded regions indicate the domain \(\Omega\),
  and the spline space spanned by the B-splines over \(\Omega\) are the same in both cases.}
  \label{fig:ghost_versus_open}
\end{figure}

In  IgA  \emph{open knot vectors} are used to
simplify the interpolation of boundary conditions, as reported  for THB in
\cite{GIANNELLI2016337} and for LRB in
\cite{johannessenIsogeometricAnalysisUsing2014}. In  \emph{open knot vectors} 
the multiplicity at  boundary knots is set to \( m = d + 1 \). An alternative approach is to
use B-splines with knot multiplicity of \( m = 1 \) along the boundary. In
order to force the partition of unity in this case, a ghost domain is added around the
domain of interest, as seen in
\cite{johannessenSimilaritiesDifferencesClassical2015a} for both THB and LRB.
This distinction is illustrated for univariate cubic splines in
\cref{fig:ghost_versus_open}.

In \cref{sec:boundary_multiplicity} we address the effects the choice of
boundary knot multiplicity has on condition numbers. To distinguish between
single multiplicity and open knots at the domain boundary we prefix any method
using single knot multiplicity on the boundary with a ghosted domain by an
``S''. Using this naming convention, the methods addressed in
\cite{johannessenIsogeometricAnalysisUsing2014} are respectively S-THB and
S-LRB. In this paper, we show that for the same tensor product spline space, THB and LRB are
superior with respect to condition numbers of mass and stiffness matrices
compared to respectively S-LRB and S-THB. We also explain the intriguing near
constant behaviour of the condition numbers reported in
\cite{johannessenSimilaritiesDifferencesClassical2015a}, where S-LRB and S-THB
were addressed, and condition numbers seemed to be nearly independent of the refinement
level. We show that this is due to single knot multiplicity at domain boundaries
for the examples presented in
\cite{johannessenSimilaritiesDifferencesClassical2015a}. Further it is showed
that for more levels of refinement the condition numbers for the mass matrix
for S-THB and S-LRB will meet and then follow the growing  curves  for
respectively THB and LRB.

In HB and THB the refinement procedure (at an element level) consists of
marking elements for splitting. Marked elements are subsequently split in both
parameter directions. This contrasts with the refinement procedure LRB allows,
namely that of splitting an element in a single parameter direction at the
time, provided that at least one B-spline is split in the process. This can be
used to modify the hierarchical refinement, and possibly improve the approximation
properties of the resulting spline space. In the remaining sections we use open
knot vectors at domain boundaries and address how such modifications influence
the condition numbers for mass and stiffness matrices for bi-cubic spline
spaces. The remaining sections are structured as follows:
\begin{description}

  \item[\cref{sec:box_partitions_meshes_spline_spaces}] gives a lightweight
    introduction to box-partitions and spline spaces over such partitions.  The
   starting point for THB and LRB refinement is a tensor product spline space.
    The key concept of element overloading is defined, the situation where more
    B-splines cover an element than are needed for spanning the polynomial
    space over the element.  We briefly summarize some key properties.
    Subsequently, we recall the definitions of both LRB and THB splines. We
    also relate the refinement strategies for LRB to 
    T-splines.

  \item[\cref{sec:local_modifications_and_overloading}] takes a look at
    overloading. We look  at  how to reduce or completely remove overloaded
    regions in a mesh. We showcase some specific overloading patterns that
    occur for hierarchical refinement of THB and LRB. Furthermore, we show how
    local modifications to the LRB-mesh reduce overloading as well as condition
    numbers.

  \item[\cref{sec:numerical_experiments}] provides a quantitative comparison
    between the methods. We conduct our numerical experiments using modified
    central and diagonal refinement examples from
    \cite{johannessenSimilaritiesDifferencesClassical2015a} with a finer
    initial tensor-product mesh. This gives enough room on each hierarchical
    level for the local modifications to take place. The examples show that LRB
    with no overloading have smaller condition numbers for the mass matrix per degree of freedom
    than THB and LRB. However, the difference seems to be so small that in
    general all methods have a similar behaviour.
  \item[ \cref{sec:conclusion}] summarizes the main results of this paper.
\end{description}

%% file: boundary_multiplicity.tex
\section{Condition numbers and knotline multiplicity at domain boundary}
\label{sec:boundary_multiplicity}

\begin{figure}[htbp]
  \centering
  \includegraphics[width=1\linewidth]{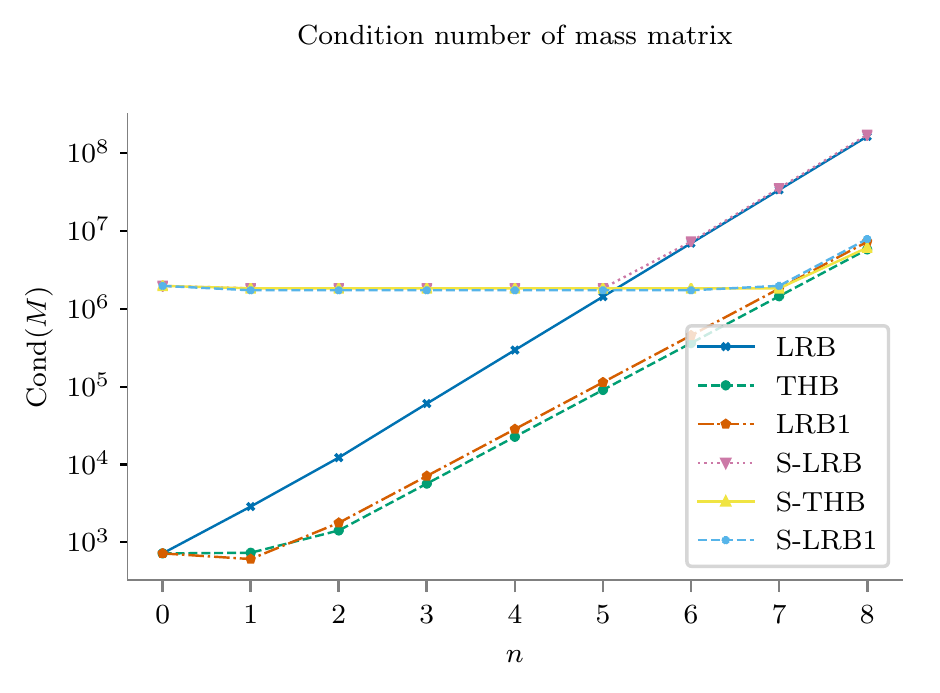}
  \caption[Condition number of mass matrix]{The condition number of the mass matrix. We see that under repeated
    refinement, the condition numbers corresponding to spline spaces with open
    knot vectors (THB, LRB) tends towards the condition numbers corresponding
    to spline spaces with single knots (S-THB, S-LRB). We also see that 
  a small local modification to reduce overloading in the LRB-space reduces the condition number of
the mass matrix (S-LRB1, LRB1).}
  \label{fig:boundary_multiplicity_mass}
\end{figure}

\begin{figure}[htbp]
  \centering
  \subfloat[]{\includegraphics[width=0.49\linewidth]{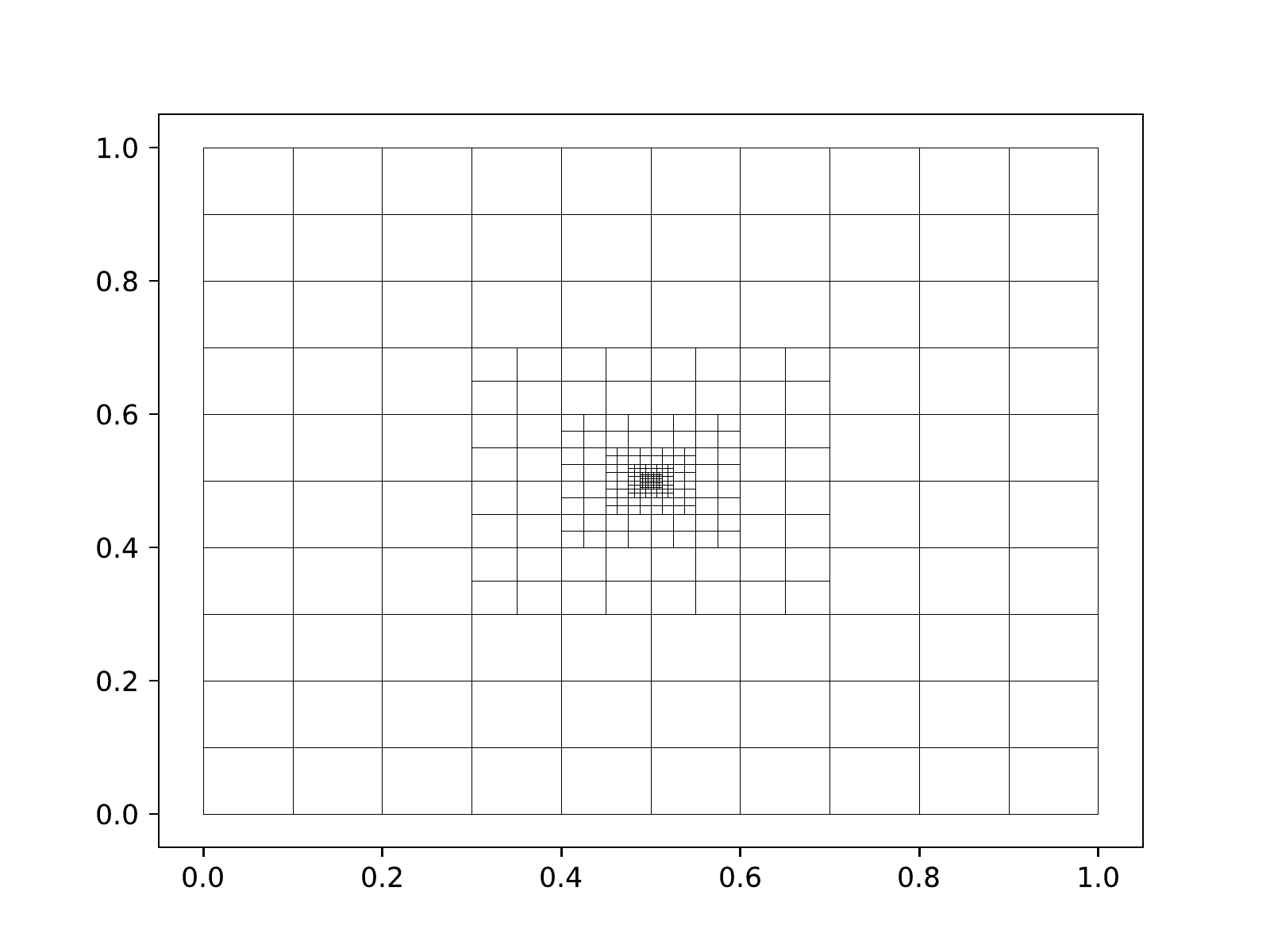}\label{sub:boundary_multiplicity_mesh_nomod}}
  \subfloat[]{\includegraphics[width=0.49\linewidth]{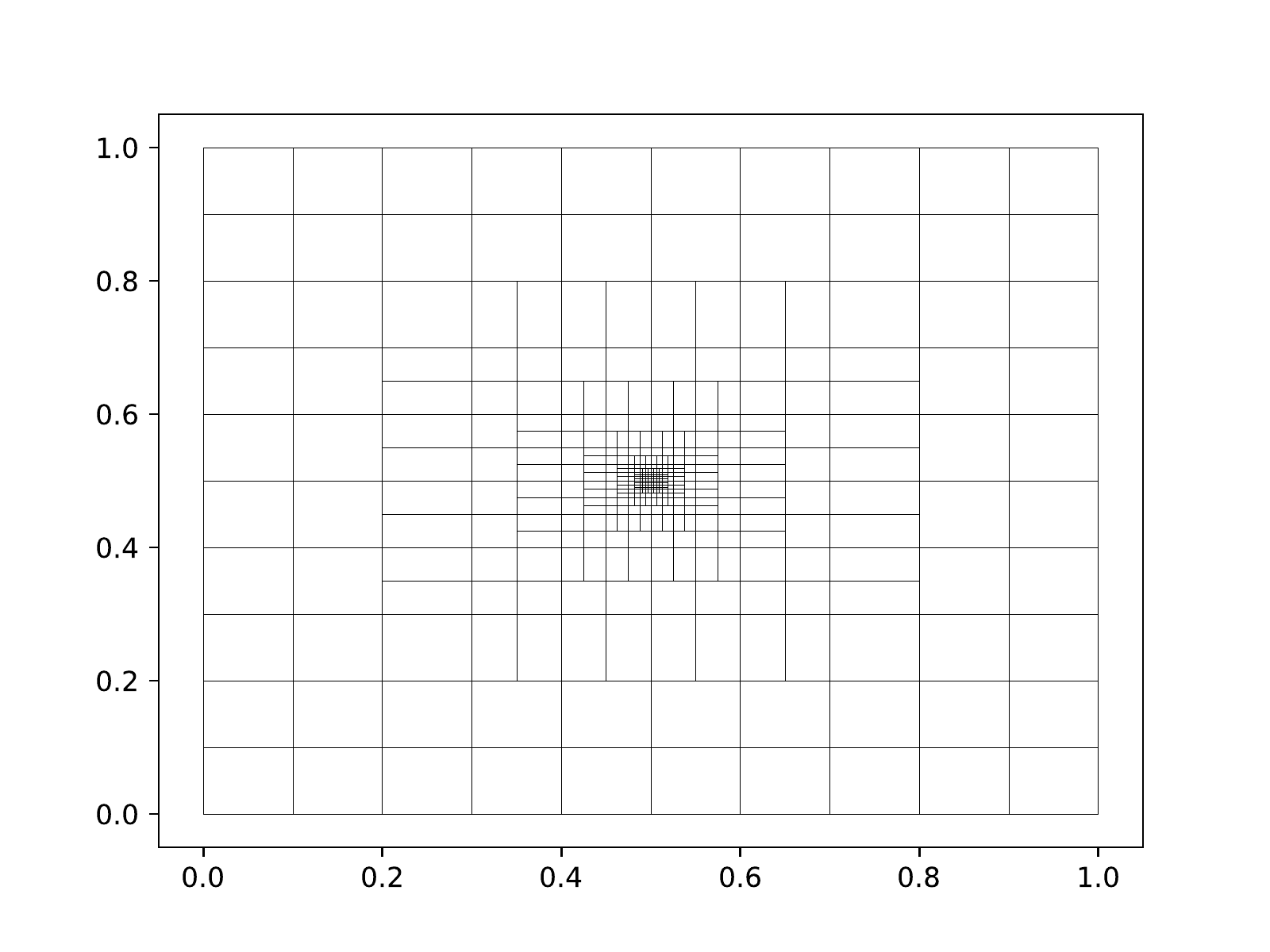}\label{sub:boundary_multiplicity_mesh_mod}
  }
  \caption[Meshes used for comparison]{The meshes used for the preliminary comparison. In
  \protect\subref{sub:boundary_multiplicity_mesh_nomod}, the unmodified mesh used for
S-THB, S-LRB, THB and LRB. In \protect\subref{sub:boundary_multiplicity_mesh_mod} the
modified mesh used for S-LRB1 and LRB1. This mesh generates a few extra degrees
of freedom.}
  \label{fig:boundary_multiplicity_mesh}
\end{figure}

In \cite{johannessenSimilaritiesDifferencesClassical2015a}, hierarchical
refinement was performed for five levels of refinement, using S-THB and S-LRB.
The results reported that the number of refinement levels had little to no
influence on the evolution of condition numbers for stiffness and mass
matrices. There were some minute differences between S-THB and S-LRB, but they
followed the same trend. In \cref{fig:boundary_multiplicity_mass} we display
the condition number of the mass matrix for up to eight refinement levels for
S-THB, S-LRB, THB and LRB when run on a hierarchical mesh from
\cite{johannessenSimilaritiesDifferencesClassical2015a}. The relevant mesh at
the fifth refinement level is displayed in
\cref{fig:boundary_multiplicity_mesh}. 

The results from \cite{johannessenSimilaritiesDifferencesClassical2015a} is
reproduced, and corresponds to the S-THB and S-LRB curves for the first five
refinements. However, at the sixth refinement, the curve for the condition
number of the mass matrix for S-LRB breaks off and grows exponentially
following the curves of LRB that starts three orders of magnitude lower.  In
\cref{fig:boundary_multiplicity_mass} there are also two additional curves
(S-LRB1 and LRB1). These are added to show that modifying the mesh by inserting
additional knot lines in one parameter direction, with the effect of reducing
overloading, significantly reduces the condition numbers of LRB-refinement.
This modified mesh is shown in \cref{sub:boundary_multiplicity_mesh_mod}. We
will discuss such modifications more closely in
\cref{sec:local_modifications_and_overloading}.

\begin{figure}[htbp]
  \centering
  \includegraphics[width=1\linewidth]{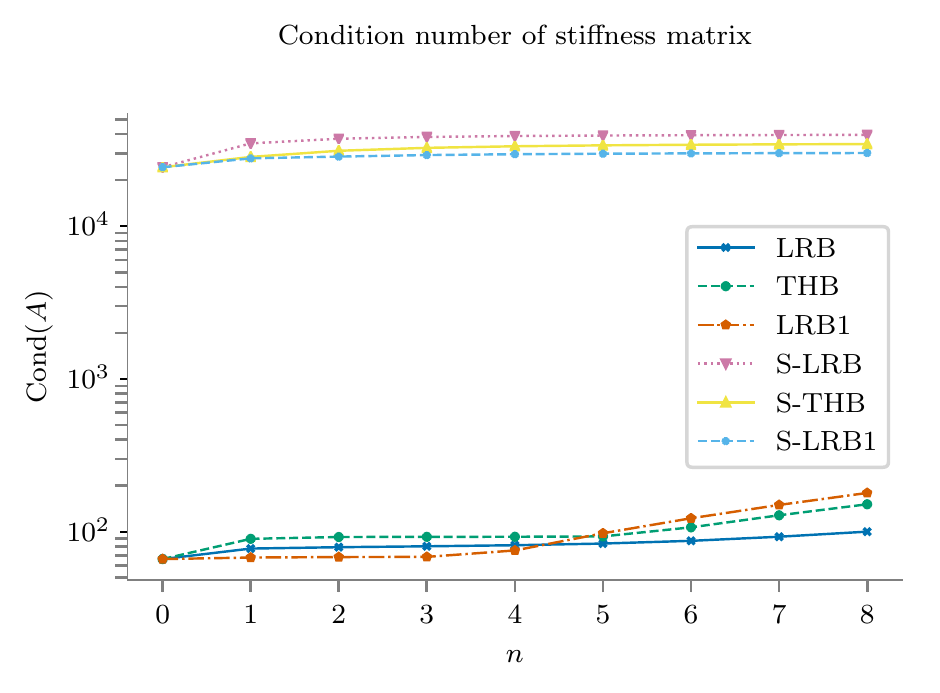}
  \caption[Condition number of stiffness matrix]{The condition number of the stiffness matrix. Here the separation
  between S-LRB, S-THB, LRB and THB are seen in even greater effect.}
  \label{fig:boundary_multiplicity_stiffness}
\end{figure}

Multiplicity of domain boundary knot lines also influences the condition number
of the stiffness matrix, as seen in \cref{fig:boundary_multiplicity_stiffness}.
Here we see that the condition numbers for single boundary knot multiplicity
(S-THB, S-LRB and S-LRB1) are two orders of magnitude higher than the condition
numbers for open knot vectors (THB, LRB, LRB-1).

\subsection{Boundary knotline multiplicities}

\begin{figure}[htpb]
  \centering
  \subfloat[\(n=1\) (LRB)]{\includegraphics[width=0.33\linewidth]{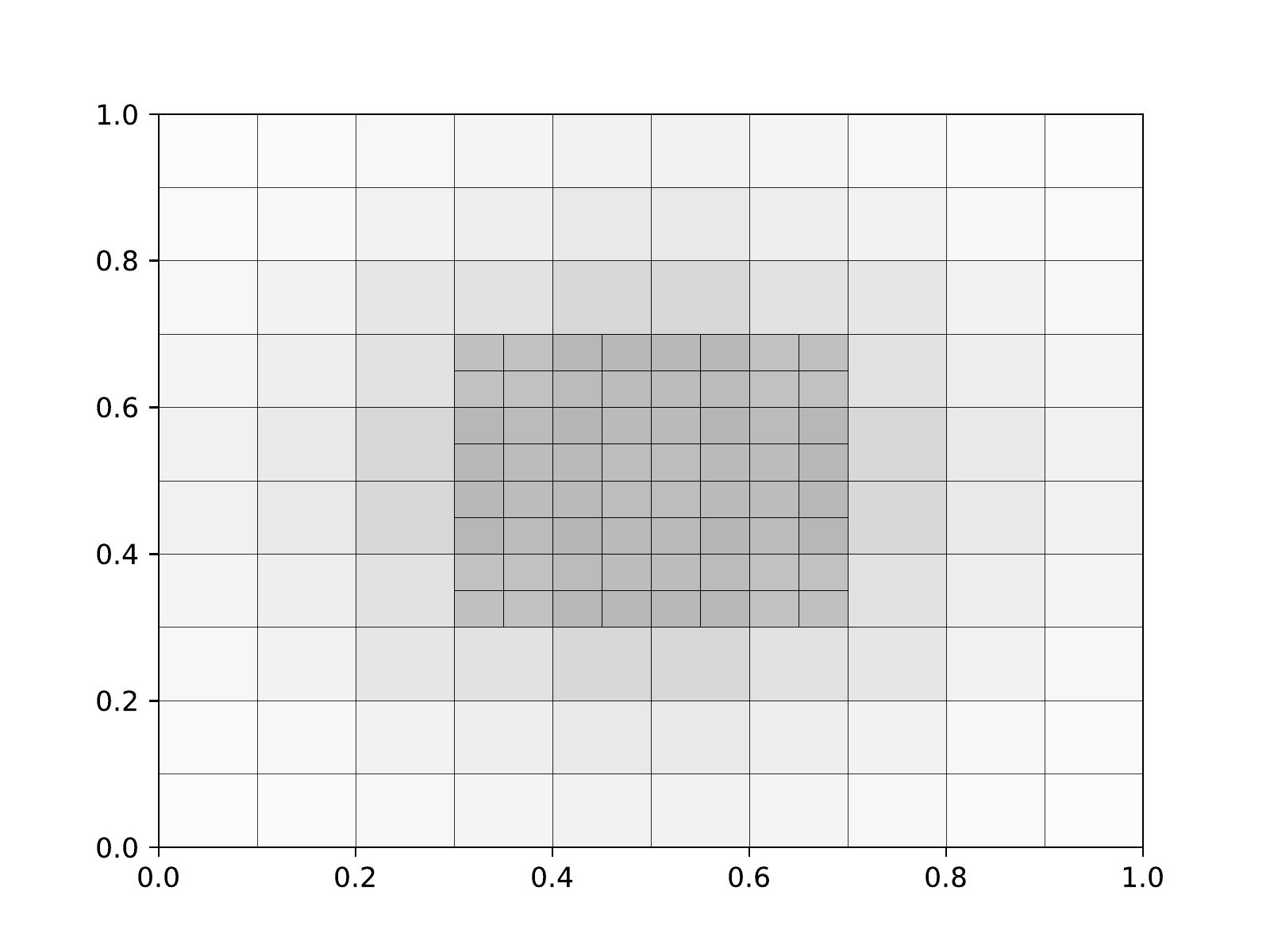}\label{sub:lrb_smallest_1}}
  \subfloat[\(n=3\) (LRB)]{\includegraphics[width=0.33\linewidth]{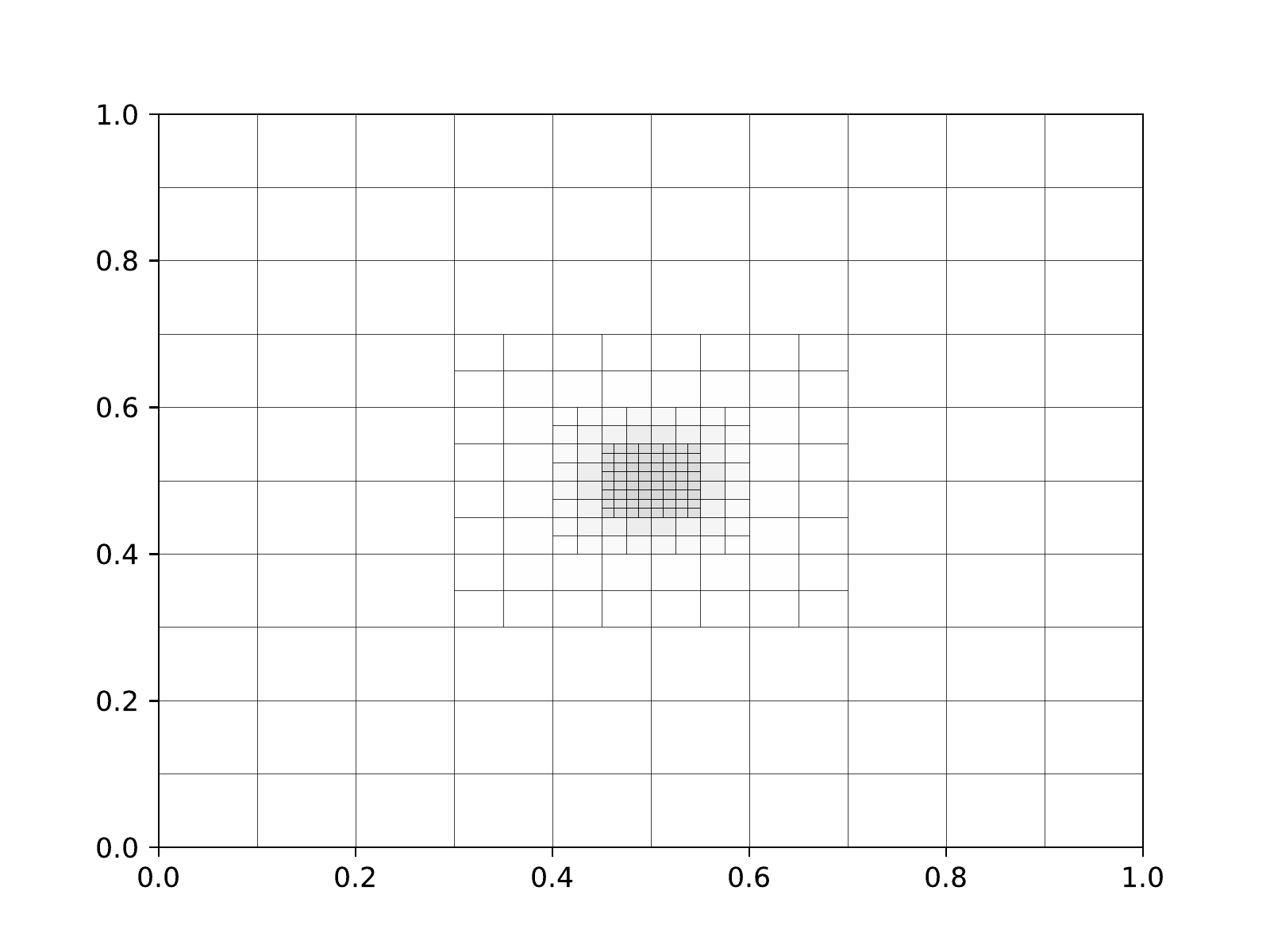}\label{sub:lrb_smallest_3}}
  \subfloat[\(n=6\) (LRB)]{\includegraphics[width=0.33\linewidth]{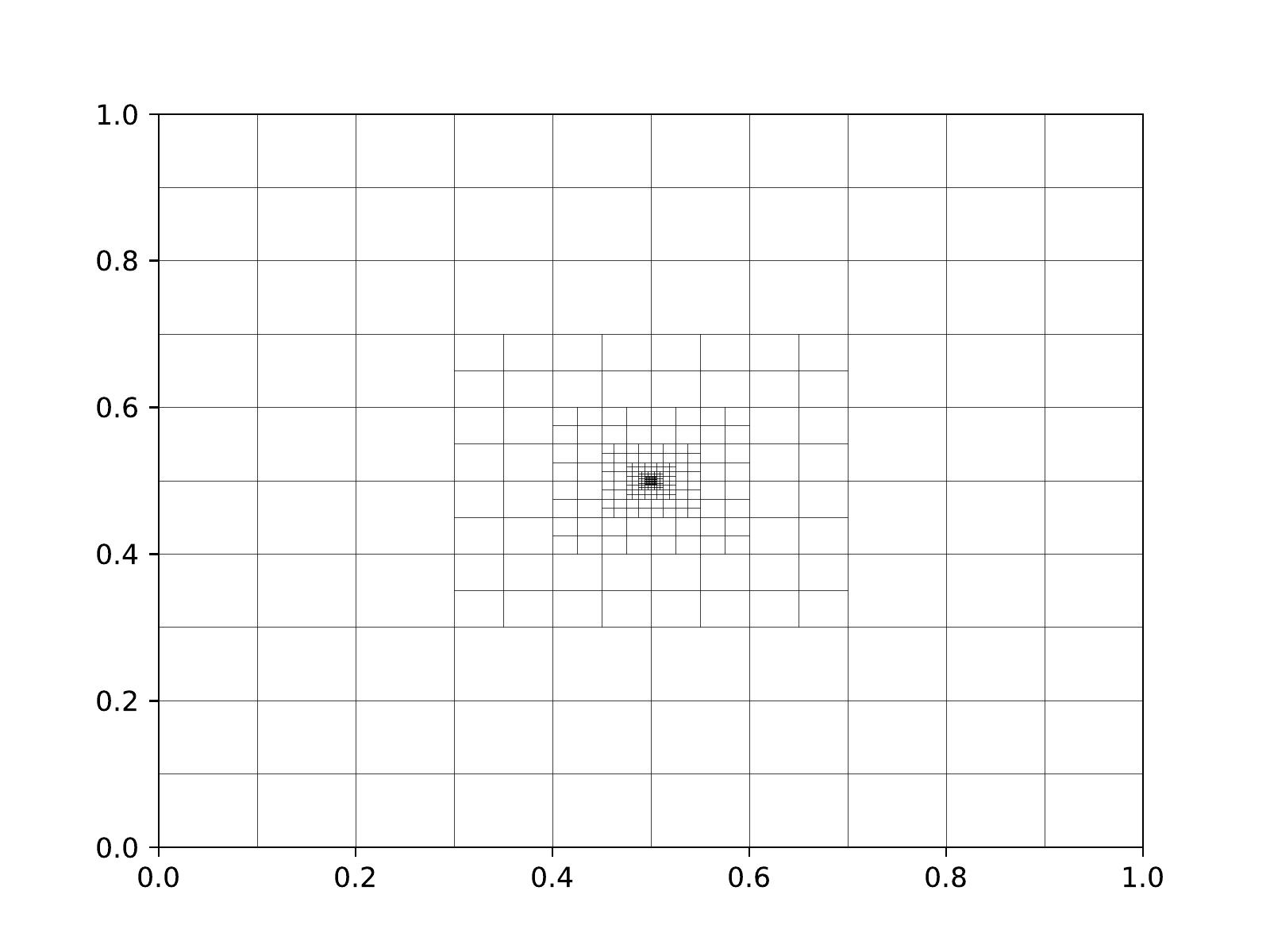}\label{sub:lrb_smallest_6}}
                   
  \subfloat[\(n=1\) (S-LRB)]{\includegraphics[width=0.33\linewidth]{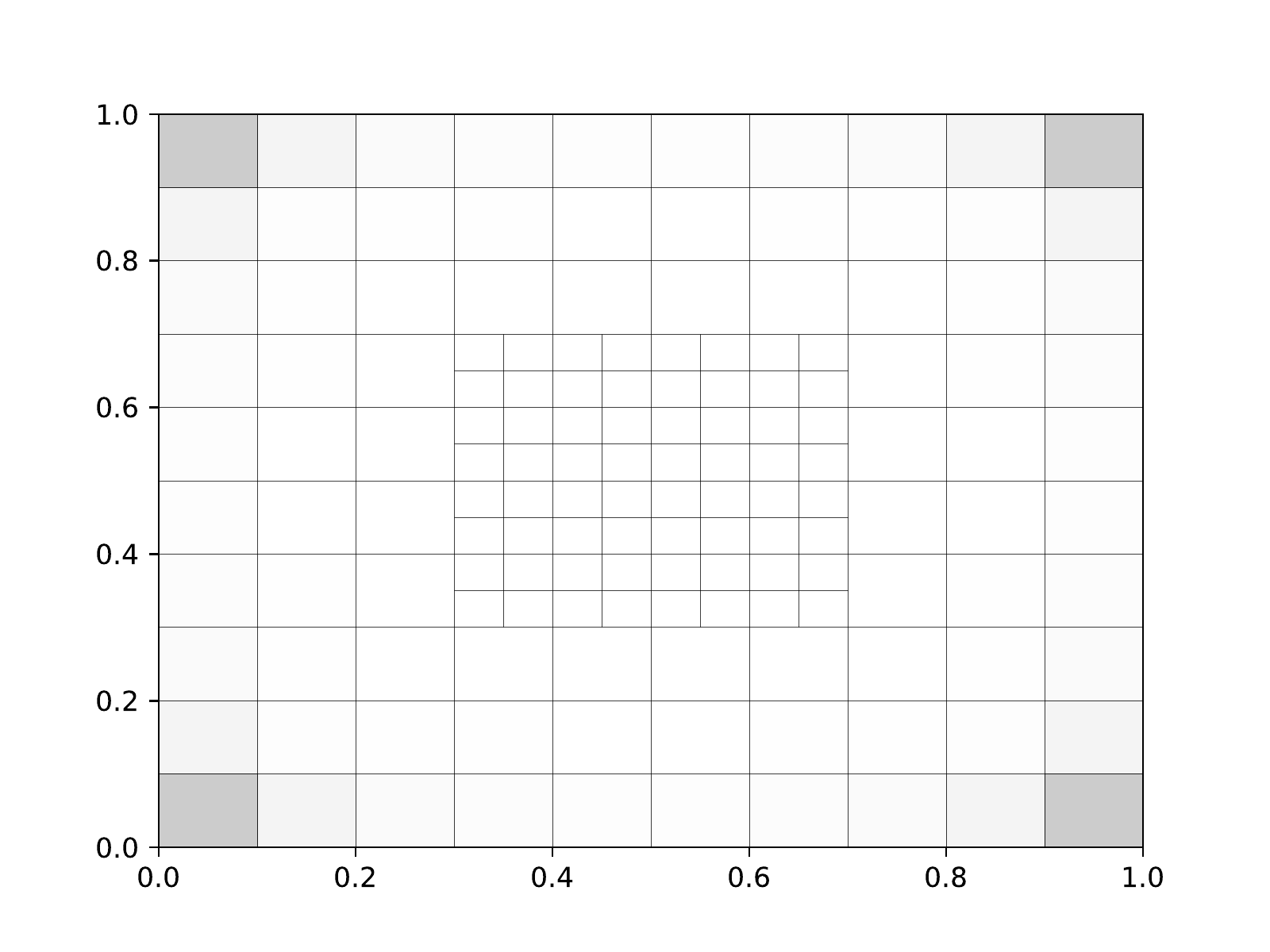}\label{sub:slrb_smallest_1}}
  \subfloat[\(n=3\) (S-LRB)]{\includegraphics[width=0.33\linewidth]{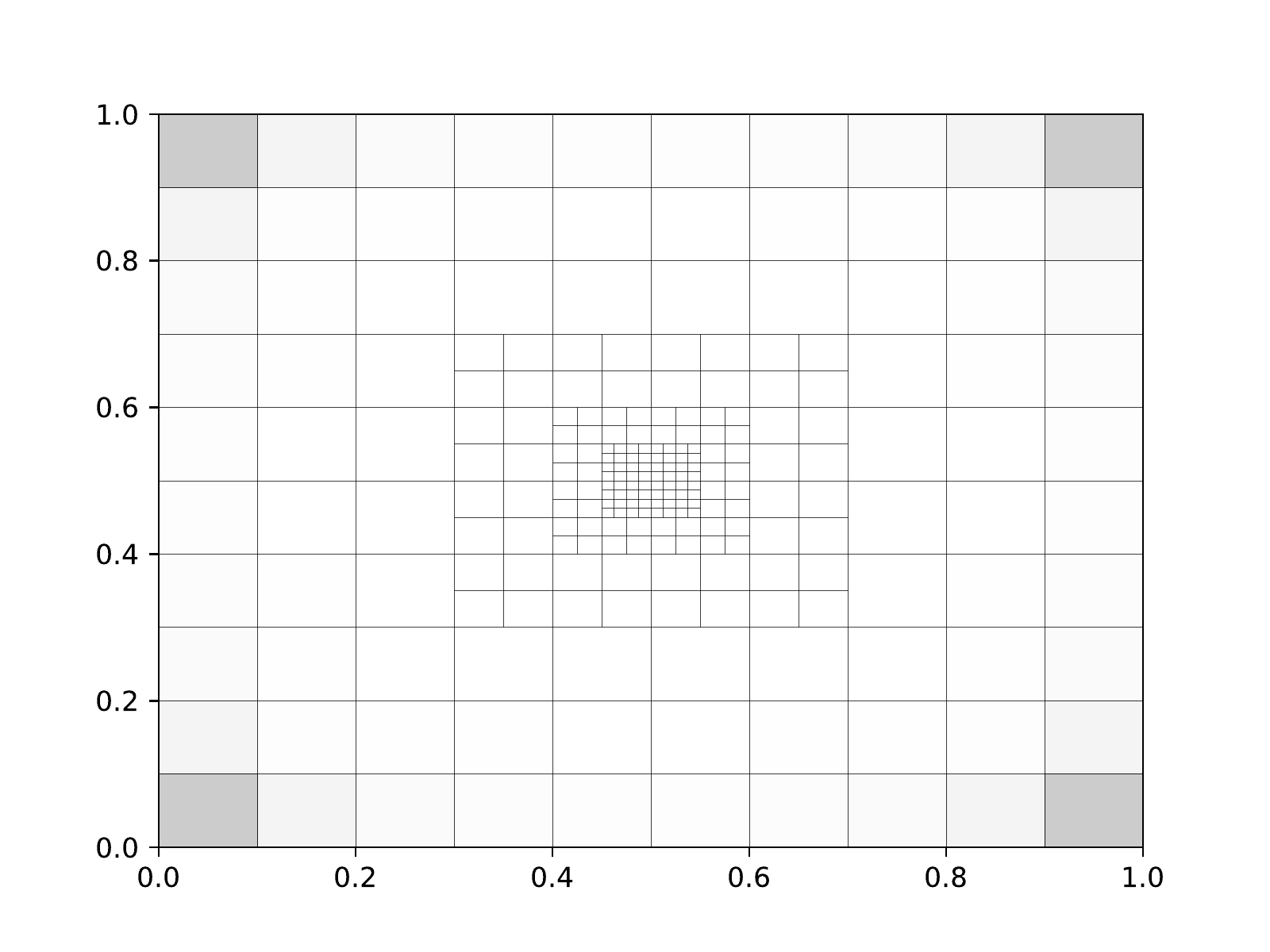}\label{sub:slrb_smallest_3}}
  \subfloat[\(n=6\) (S-LRB)]{\includegraphics[width=0.33\linewidth]{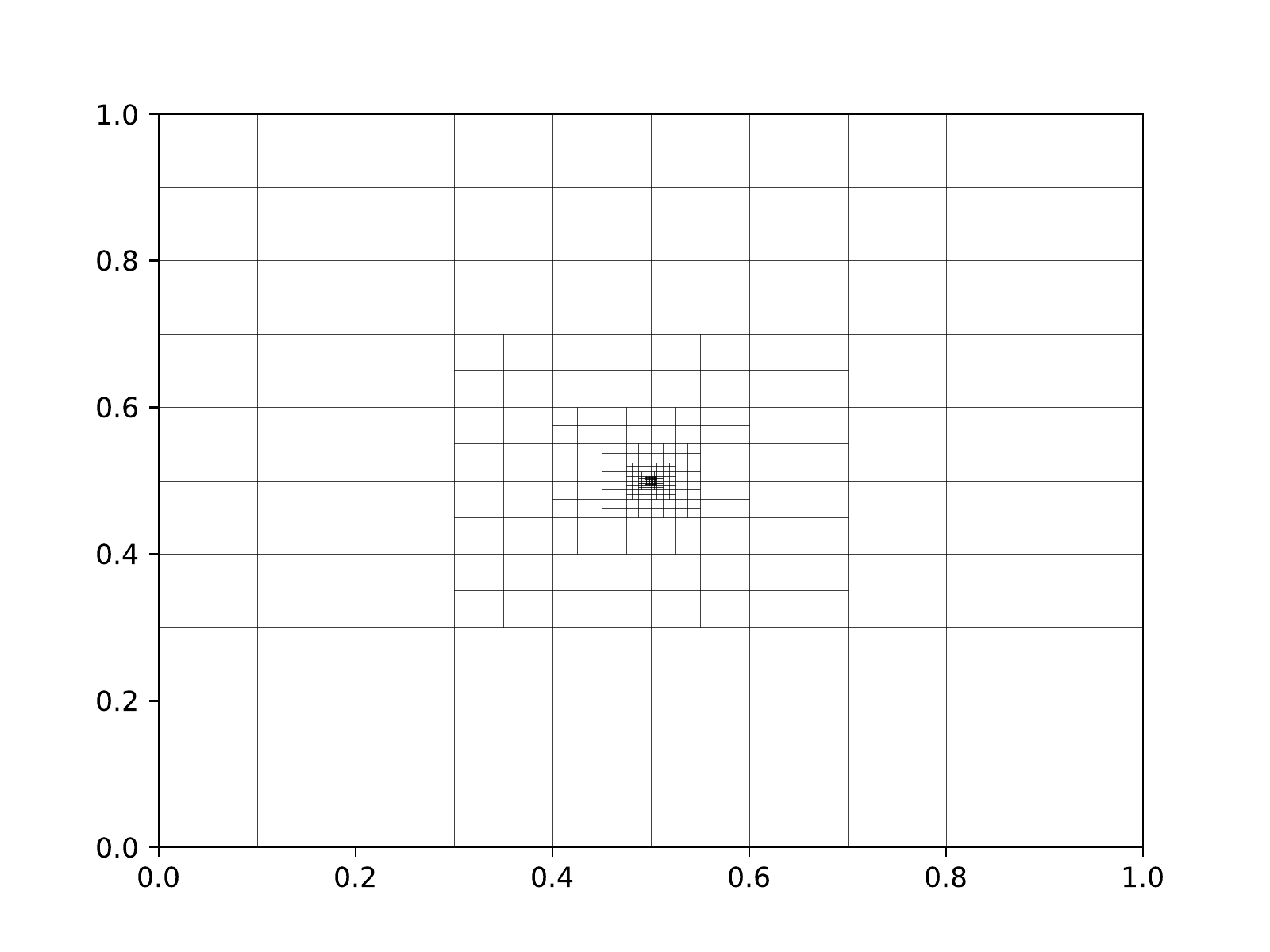}\label{sub:slrb_smallest_6}}
  \caption[Eigenvectors corresponding to the smallest eigenvalues]{The
    eigenvectors corresponding to the smallest eigenvalue of the mass matrix
    for LRB (top row), and for S-LRB (bottom row) visualized over the
    hierarchical mesh after one, three and six refinements (left to right).
    Darker color indicates higher influence. As we see, the smallest
    eigenvalues for LRB is localized in the refined region after only one
    refinement. On the other hand, S-LRB is localized in the corners of the
    domain up until but not including six refinements, as shown for \( n = 1 \)
    and \( n = 3 \). The effect of the locally refined region dominates only
  after \( n = 6 \) refinements as in \protect\subref{sub:slrb_smallest_6}.} 
  \label{fig:smallest_eigenvector_mass}
\end{figure}

\begin{figure}[htpb]
  \centering
  \subfloat[\(n=1\) (LRB)]{\includegraphics[width=0.33\linewidth]{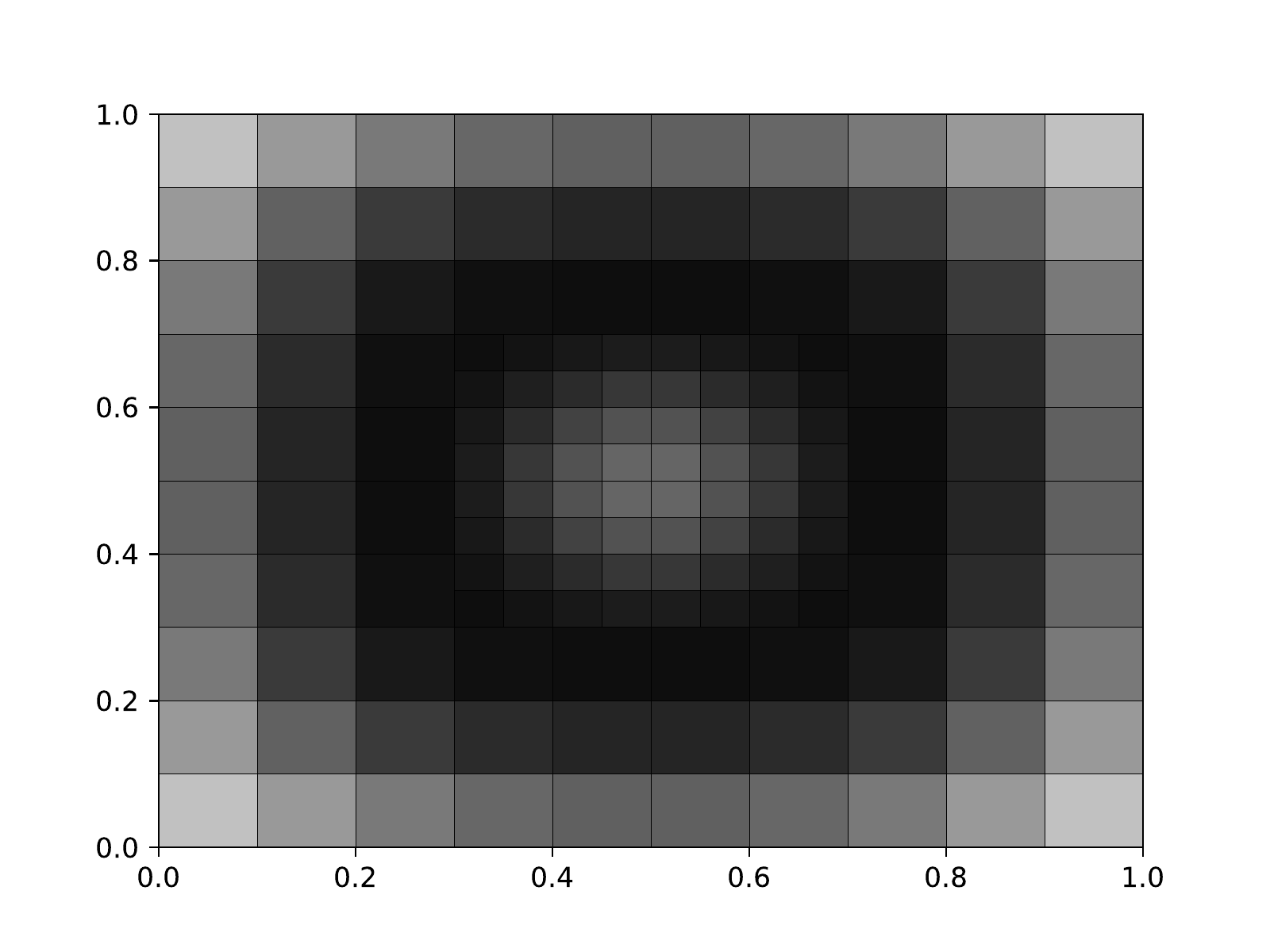}\label{sub:lrb_largest_1}}
  \subfloat[\(n=3\) (LRB)]{\includegraphics[width=0.33\linewidth]{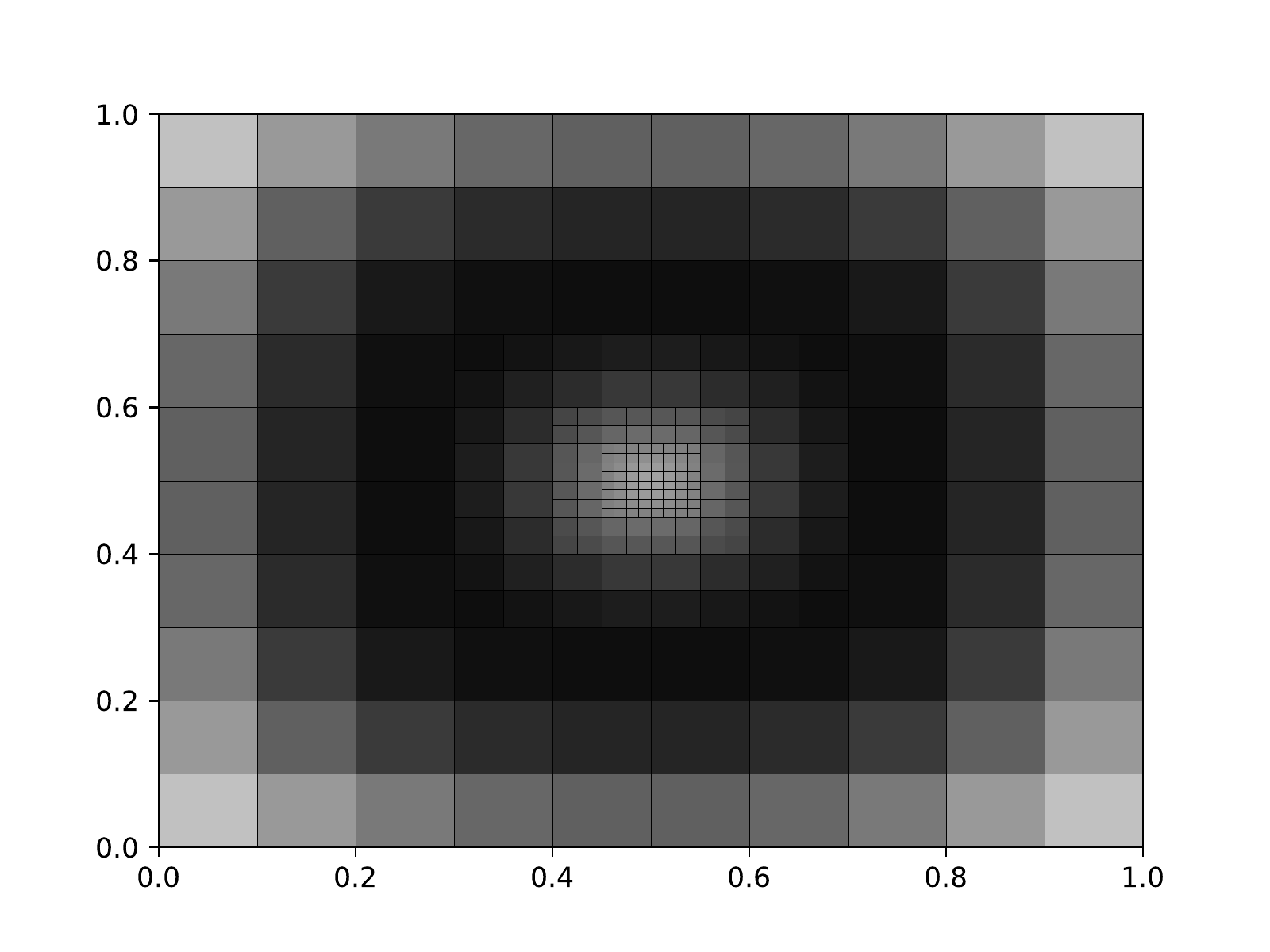}\label{sub:lrb_largest_3}}
  \subfloat[\(n=6\) (LRB)]{\includegraphics[width=0.33\linewidth]{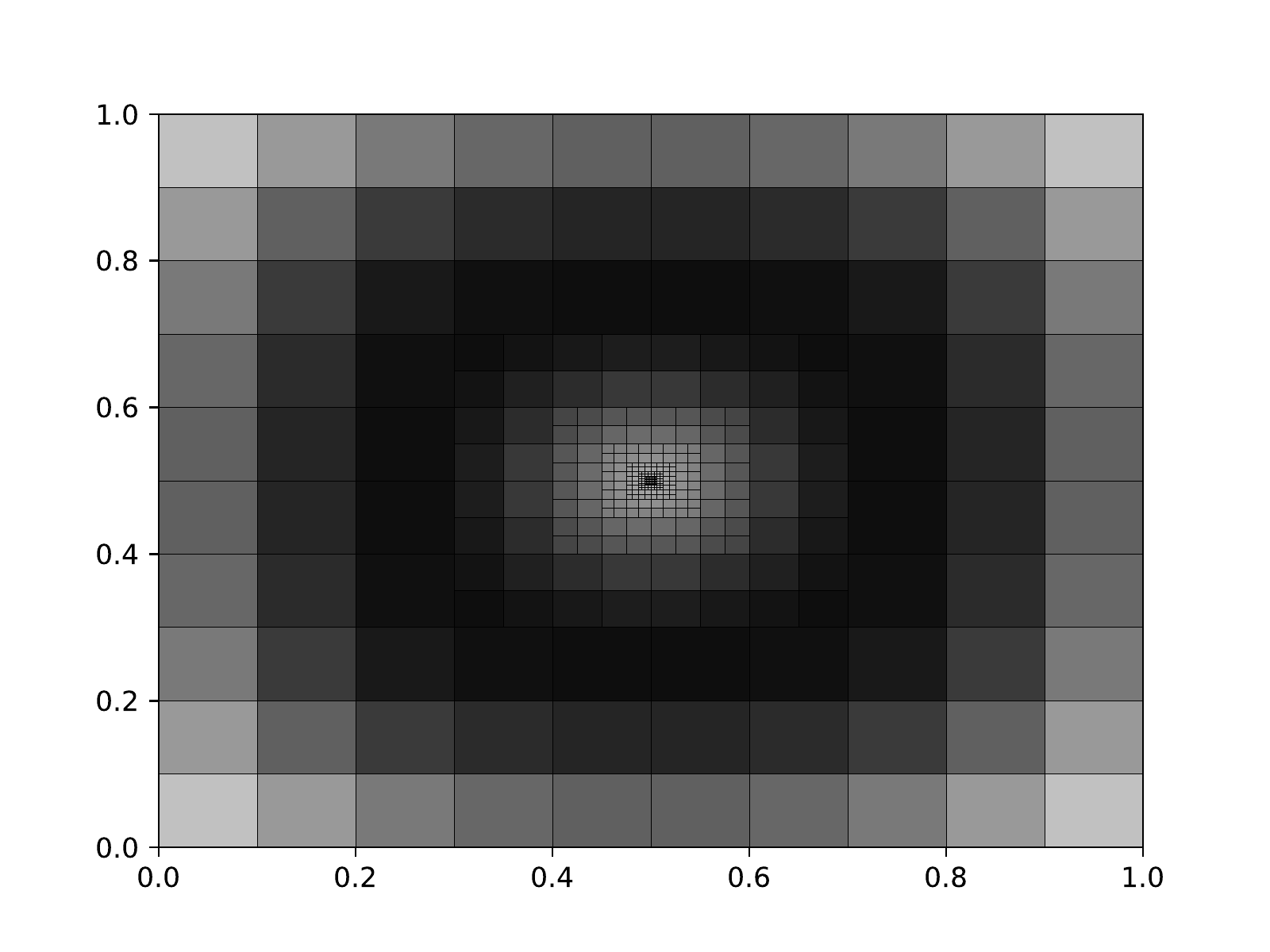}\label{sub:lrb_largest_6}}

  \subfloat[\(n=1\) (S-LRB)]{\includegraphics[width=0.33\linewidth]{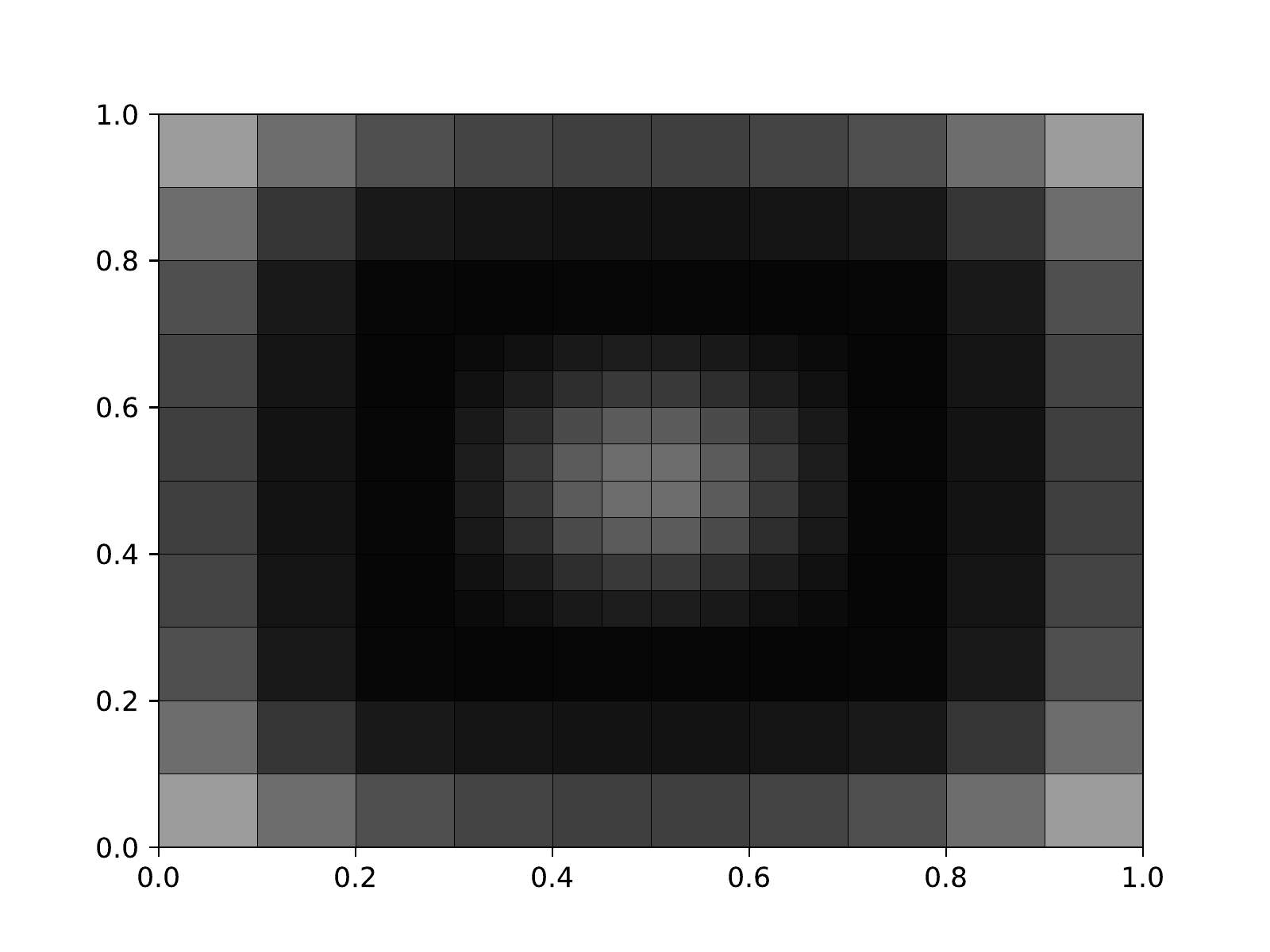}\label{sub:slrb_largest_1}}
  \subfloat[\(n=3\) (S-LRB)]{\includegraphics[width=0.33\linewidth]{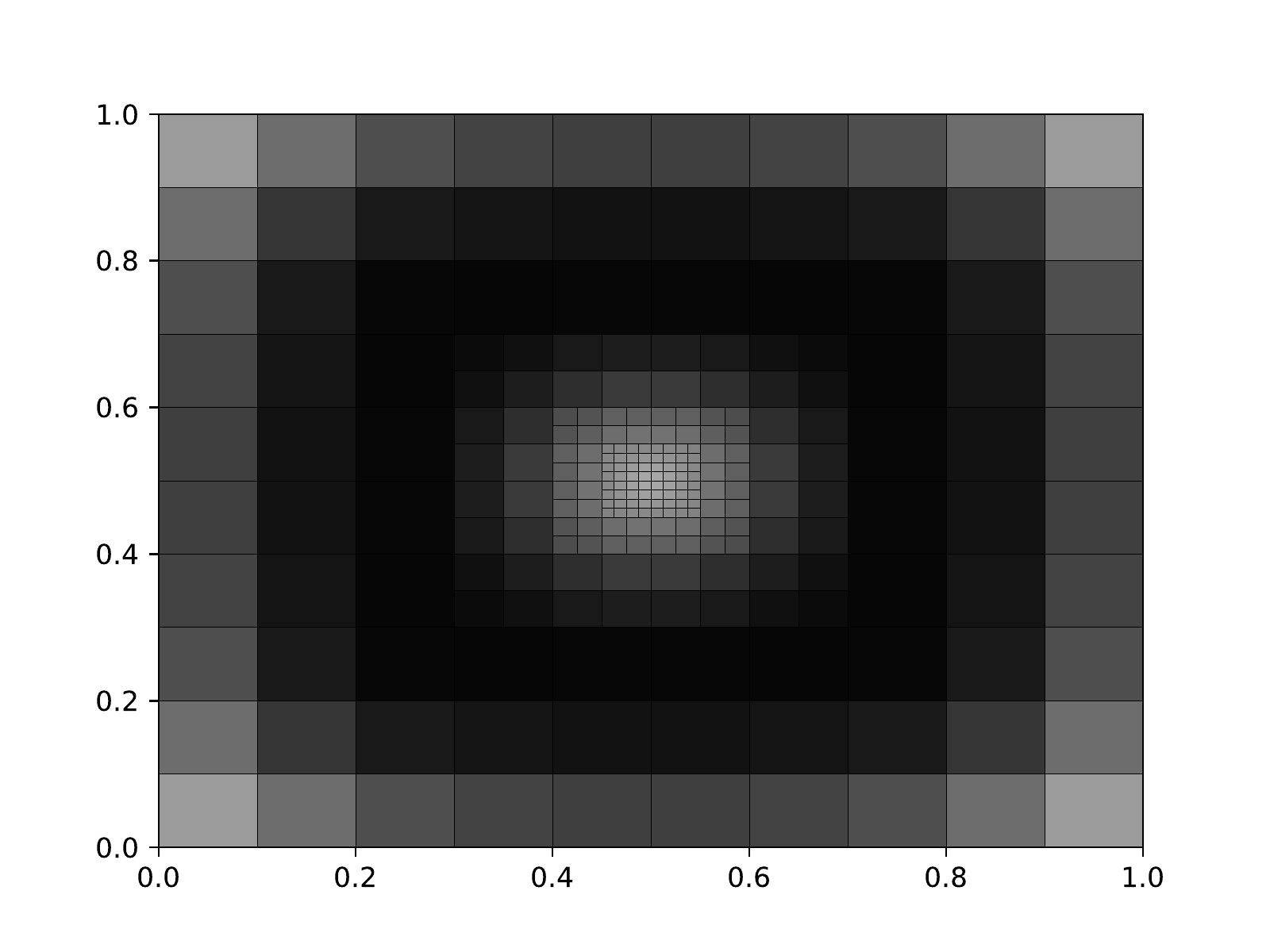}\label{sub:slrb_largest_3}}
  \subfloat[\(n=6\) (S-LRB)]{\includegraphics[width=0.33\linewidth]{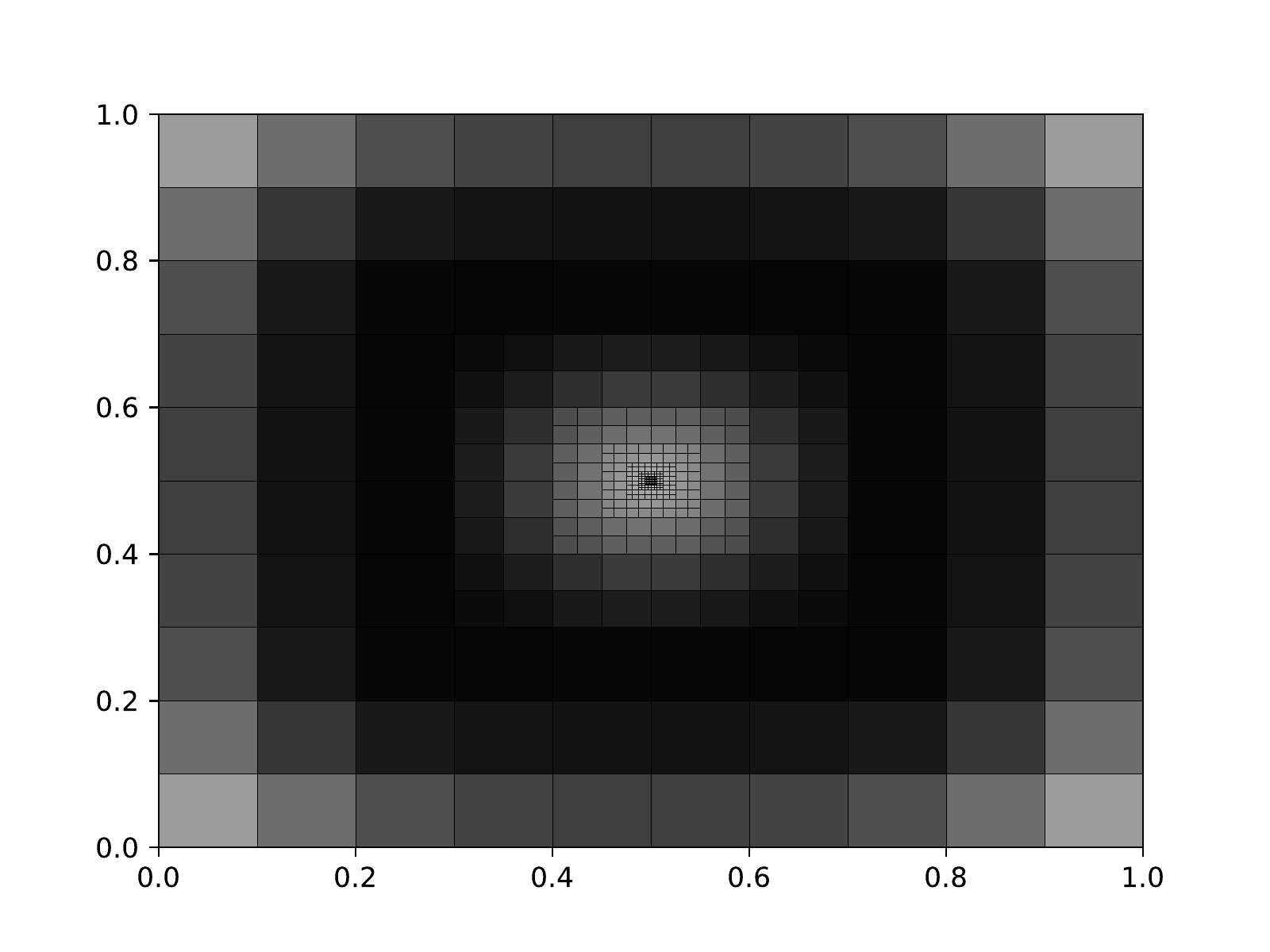}\label{sub:slrb_largest_6}}
  \caption[Eigenvectors corresponding to the largest eigenvalues]{ The
    eigenvectors corresponding to the largest eigenvalue of the mass matrix for
    LRB (top row) and for S-LRB (bottom row) visualized over the hierarchical
    mesh after one, three and six refinements (left to right).  Darker color
  indicates higher influence.  For the largest eigenvalue, the two methods
approximately correspond geometrically, and the largest eigenvalues are
constant over refinement levels for each of the methods.} 
  \label{fig:largest_eigenvector_mass}
\end{figure}

We now take a stab at explaining the drastic change in behaviour occuring at \(
n = 6 \) refinements for the S-LRB and S-THB methods as shown in
\cref{fig:boundary_multiplicity_mass}. Since the condition number of a matrix
are computed in terms of its largest and smallest eigenvalues, we decided to
take a look at the geometric localization of the eigenvectors corresponding to
these eigenvalues. By coloring the hierarchical mesh based on the influence of
each in terms of the corresponding coefficient in the eigenvector, we obtained
a rudimentary geometric visualization of these eigenvectors. In
\cref{fig:smallest_eigenvector_mass}, we see the smallest eigenvector for the
mass-matrix corresponding to  LRB and S-LRB at the first, third and sixth
refinement, and in \cref{fig:largest_eigenvector_mass}, the corresponding
largest eigenvector.

These figures correspond to the behaviour observed in
\cref{fig:boundary_multiplicity_mass} where the conditioning for LRB grows
after only one refinement, whereas S-LRB needs six refinements before the
behaviour in the refined region is registered.

The reason for this behaviour is due to the size of the B-splines defined along
the boundary in comparison to the size of the B-splines defined in the interior
of the domain. In order to illustrate this, we compute analytically the entries
in the mass matrix corresponding to B-splines on various tensor product level
and compare these values to the mass matrix entry corresponding to a B-spline
defined in the corner of the domain. 

\subsubsection*{Observation for the mass matrix}

Over the domain \( \Omega = [0, 1] \times [0, 1] \) we define a tensor product
grid with element size \( \ell \). In the case of bi-cubic spline spaces, the
B-spline defined in the lower left corner of the domain can for LRB and S-LRB be written in terms of their knots as
\begin{align}
  \begin{split}
    B &\coloneqq B[\vec{x}] B[\vec{y}], \\
    Q &\coloneqq B[\vec{s}] B[\vec{t}],
  \end{split}
\end{align}
where \( \vec{x} = \vec{y} = [0, 0, 0, 0, \ell] \) and \( \vec{s} = \vec{t} =
[-3\ell, -2\ell, -\ell, 0, \ell] \). In both cases, the two B-splines have only
one element of support in the domain \( \Omega \), namely \( \element \coloneqq [0,
\ell] \times [0, \ell] \). To get a feel for the differences in
influence on the mass matrix these B-splines have, we compute the corresponding diagonal
elements in the mass matrix.

The polynomial restrictions of \( B \) and \( Q \) to the element \( \element \) is
\begin{align}
  \begin{split}
    \eval[1]{B}_{\element} (x, y) &= \frac{(\ell - x)^3(\ell - y)^3}{\ell^6}, \\ 
    \eval[1]{Q}_{\element} (x, y) &= \frac{(\ell - x)^3(\ell - y)^3}{36\ell^6}\,.
  \end{split}
\end{align}
In other words, \( \eval[0]{B}_{\element} = 36 \eval[0]{Q}_{\element} \). If we now compute
the diagonal mass matrix entries corresponding to these two elements, we obtain
the following:
\begin{equation}
    \int_\element B^2 = \frac{\ell^2}{49}, \quad \int_\element Q^2 = \frac{\ell^2}{49} \cdot \frac{1}{36^2}\, .
\end{equation}
We here see that the matrix element corresponding to the corner B-spline \( Q
\) arising in S-LRB is three orders of magnitude smaller than the matrix
element corresponding to B. Recall the disparity between the curves in
\cref{fig:boundary_multiplicity_mass}, where the differences in the condition
numbers also were three orders of magnitude.

\begin{remark}
  This relation between \( Q \) and \( B \) is both dimension and
  degree-dependent. The effect will be magnified for higher spatial dimension
  and higher polynomial degree. 
\end{remark}

%% file: box_partitions_meshes_spline_spaces.tex
\section{Box partitions, meshes, and spline spaces}
\label{sec:box_partitions_meshes_spline_spaces}

In this section we review the concept of \emph{box partitions} and spline
spaces over such partitions. While the construction generalizes to any
dimension, we will gradually focus our attention to the two-dimensional case,
as this is most relevant for our discussion. A fully general treatment can be
found in
\cite{dokkenPolynomialSplinesLocally2013,patriziLinearDependenceBivariate2018}.
The fundamental building block of a box partition is the
\emph{\(d\)-dimensional box}.

\begin{definition}
A \define{box} \(\dbox\) in \(\R^d\) (or \define{\(d\)-box}) is the Cartesian product of \( d \) closed
finite intervals \( J_1, \ldots, J_d\):
\begin{equation}
\dbox = \bigtimes_{i=1}^{d} J_i.
\end{equation}
The \define{dimension} of \( \dbox \) is defined to be the number of
non-trivial intervals in its definition, and is denoted \( \dim(\dbox) \). We
call a \(d\)-box of dimension \(d\) an \define{element}, while a \( d\)-box
of dimension \( d-1 \) is called a \define{mesh-rectangle}. To any
mesh-rectangle, we associate an integer \( k \) corresponding to which
parametric dimension its trivial component resides, and we call the
mesh-rectangle a \define{\(k\)-mesh-rectangle} if this has to be emphasized.
\end{definition}

In the two-dimensional setting \((d = 2)\), a \define{meshline} is a
one-dimensional mesh-rectangle.
\begin{remark}
Note that these naming-conventions are independent of the dimension of the
ambient space. Hence, a \emph{mesh-rectangle} may very well be something
different from a rectangle. As an example, a mesh-rectangle in \( \R^4 \) is
an axis aligned box. Furthermore, the integer \( k \) corresponding to any
mesh-rectangle encodes the \emph{direction} of the mesh-rectangle. In the
two-dimensional case, where mesh-rectangles are lines, a \(1\)-mesh-rectangle
is a vertical line, and a \(2\)-mesh-rectangle is a horizontal line.
\end{remark}

As customary in discretizations of computational domains, a large domain is
partitioned into a set of non-overlapping smaller geometrical entities. We call
such a partition in this specific setting a \emph{box partition}, and this is
more precisely defined as follows:

\begin{definition}
Let \( \Omega \subset \R^d \) be an element (\(d\)-box of dimension \(d\)).
A finite collection \(\boxpartition\) of elements is said to be a \define{box
partition} of \( \Omega \) if
\begin{enumerate}
\item\(\interior{\element}_i \cap \interior{\element}_j = \emptyset \text{ for all } \element_i,
\element_j \in \boxpartition \text{ where } \element_i \neq \element_j \), and
\item\(\displaystyle\bigcup_{\element \in \boxpartition} \element = \Omega \,\).
\end{enumerate}
In other words, a box partition is an interior-disjoint partition of \(
\Omega \) into a set of smaller elements.
\end{definition}

Associated to any element \( \element \) is its \emph{boundary}, which
naturally consists of boxes of dimension one less, i.e., mesh-rectangles.
Given a box partition of a larger element \( \Omega \), it is therefore
sensible to discuss the set of mesh-rectangles associated to this box
partition.

\begin{definition}[Informal]
Given a box partition \( \boxpartition \) of a domain \( \Omega \), we may
naturally associate a set of mesh-rectangles \( \mesh \) called a \define{box
mesh} on \( \Omega \) formed by taking unions of element boundaries.
\end{definition}

\begin{remark}
The link between a box partition \( \boxpartition \) and the associated box
mesh \( \mesh \) is such that by knowing one of them you may recover the
other. The box mesh \emph{generated} by a box partition is denoted \(
\mesh(\boxpartition) \), and the box partition \emph{generated} by a box mesh
is denoted \( \boxpartition(\mesh) \).
\end{remark}

As our ultimate goal is to define spline spaces based on tensor-product splines
over box-partitions, we need to have a concept of knot multiplicity in this
more general setting. 

\begin{definition}
A \define{box mesh with multiplicity} is a pair \( (\mesh, \mu) \) where \(
\mu\colon \mesh \to \N \) associates to each mesh-rectangle \( \meshline \) a
positive integer \( \mu(\meshline) \), called the \define{multiplicity} of
the mesh-rectangle. Note that this is completely analogous to the notion of
knot multiplicity for univariate B-splines.
\end{definition}

\begin{definition}
Let a polynomial multi-degree \( p = (p_1, \ldots, p_d) \) as well as a box
mesh with multiplicity \((\mesh, \mu)\) corresponding to the box-partition \(
\boxpartition \) of a \(d\)-dimensional domain \( \Omega \) be given. The
\define{spline space of degree \(p\) over \( \mesh \)} is defined as
\begin{multline}
\splinespace{\mesh}{p}{\mu} \coloneqq \Big\{f \colon \Omega \to \R :
\eval[0]{f}_{\element} \in \polynomials_p \text{ for all } \element \in 
\boxpartition \\ \text{ and } f \in \continuous{p_k - \mu(\meshline)} \text{ for
all \( k\)-mesh-rectangles } \meshline \in \mesh, \\ \text{ with }k = 1, \ldots, d\Big\}.
\end{multline}
\end{definition}

A dimension formula for general spline spaces over box partitions was presented
in \cite{pettersenDimensionMultivariateSpline2013}. In general, the dimension
depends on both the topological properties of the box partition and the
parametrization of the box partition. In the two-dimensional case --- with some
requirements on the length of the constituent meshlines --- the formula reduces
to a formula depending only on the topological features of the mesh. We
consider this outside the scope of this text, and refer the reader to
\cite{pettersenDimensionMultivariateSpline2013} for details.

In order to compute with spline spaces over box partitions of the form above,
we must be able to construct a set of basis functions that span this space.
Several constructions has been studied. We will only be dealing with 
Truncated Hierarchical B-splines, and Locally Refined B-splines.

Before we move on, we define the notion of a \emph{hierarchical mesh}, a type of
box partition over which spline bases of the aforementioned type may be
defined. This provides a common ground for comparison of the two methods. The
construction is simple and relies on marking regions for which a tensor product
mesh of various refinement levels is used.

\begin{figure}[htpb]
\centering
\subfloat[]{\includegraphics[width=0.49\linewidth]{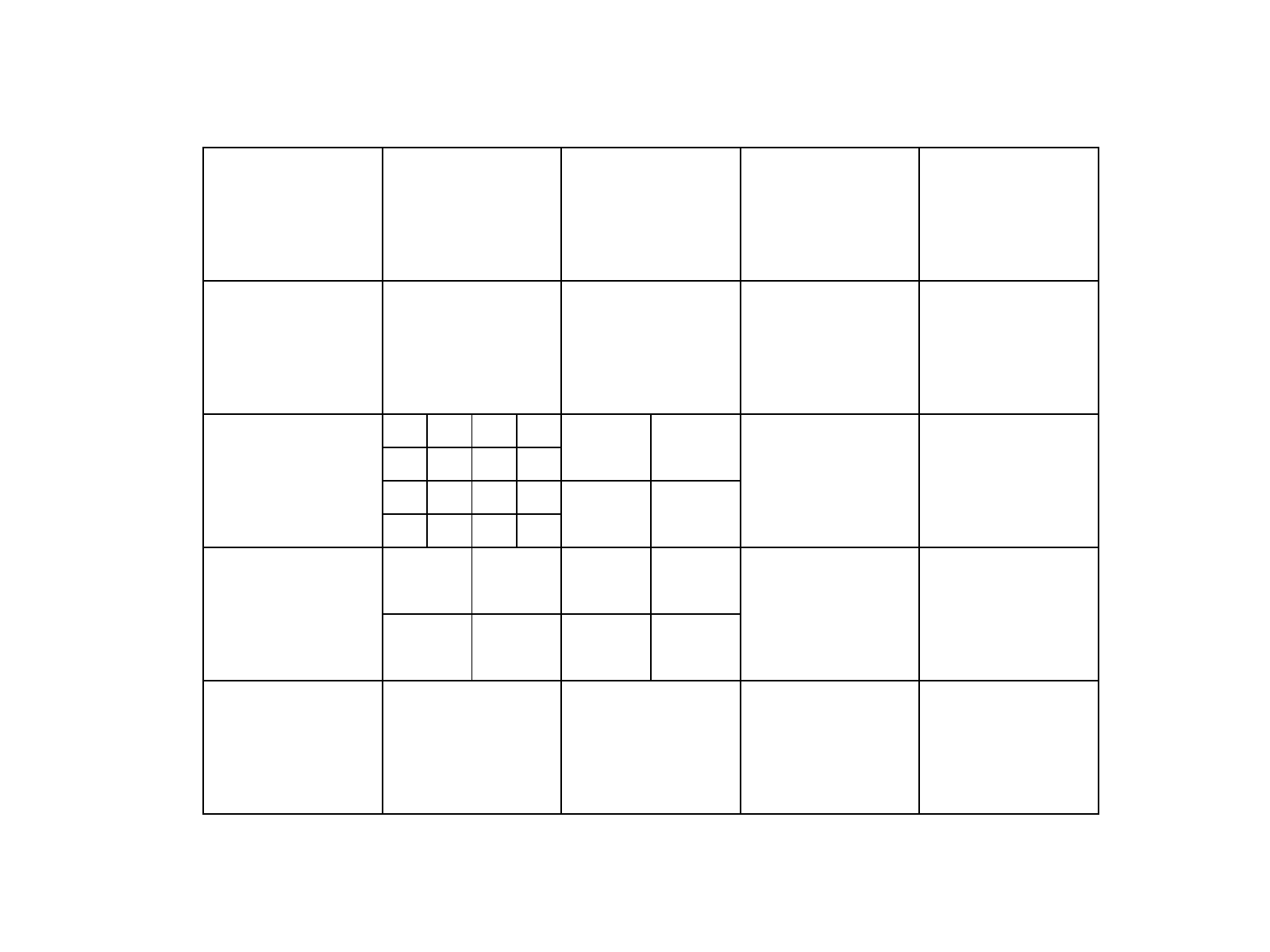}\label{sub:hierarchical_example_a}}
\subfloat[]{\includegraphics[width=0.49\linewidth]{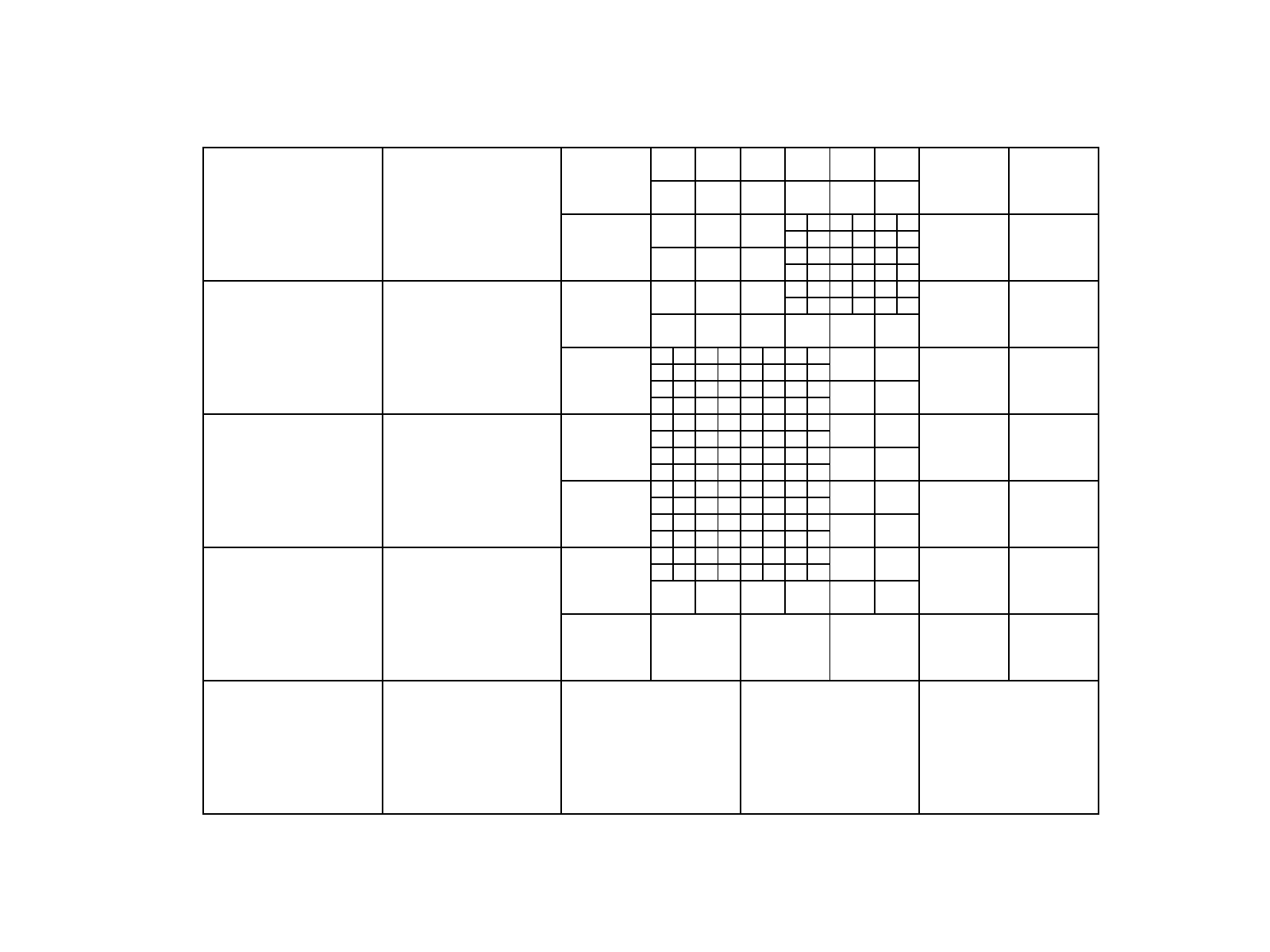}\label{sub:hierarchical_example_b}}
\caption[Examples of hierarchical meshes]{Two examples of hierarchical meshes. In \protect\subref{sub:hierarchical_example_a} a mesh
consisting of three levels of refinement, and in \protect\subref{sub:hierarchical_example_b} a mesh with
four levels of refinement. Note here that the region residing at level \( \ell
= 4 \) consists of two disjoint components.}
\label{fig:hierarchical_example}
\end{figure}

\begin{definition}
\label{def:hierarchical_mesh}
Let \( \Omega \) be a domain, and let \( \mesh_1 \subset \cdots \subset \mesh_M
\) be a sequence of nested tensor product meshes on \( \Omega \). Let \(
\Omega_1 \supset \cdots \supset \Omega_M \) be a set of nested subsets of \(
\Omega \) whose boundaries \( \boundary\Omega_\ell \) align with the meshlines
of the corresponding mesh on the coarser level \( \mesh_{\ell - 1} \) for \(
\ell = 2, \ldots, M \). The \define{hierarchical mesh} \( \mesh \) is defined as
\begin{equation}
\mesh = \set{\meshline \cap \Omega_\ell : \meshline \in \mesh_\ell \text{ for } \ell = 1, \ldots, M}, 
\end{equation}
i.e., \( \mesh \) consists of meshlines from each level intersected with the
corresponding region, see \cref{fig:hierarchical_example}.
\end{definition}

\subsection{Tensor product splines}

The foundation for all the locally refined spline spaces over box partitions
addressed in this paper, is the tensor product B-spline. We start by glossing
over some preliminary definitions.

Recall that a univariate B-spline of polynomial degree \( p \) relies on
exactly \(p + 2\) knots. This observation enables us to define B-splines
locally without referring to some global knot vector.

\begin{definition}
Given a polynomial degree \( p \) and a non-decreasing knot-vector \(
\knotvector = ( \knot_1, \ldots, \knot_{p+2}) \), we recursively define the
\define{univariate B-spline} \( B[\knotvector]\colon \R \to \R \) as follows: 
\paragraph{If \( p = 0 \), then} 
\begin{equation}
B[\knotvector] = \begin{cases}
1, & x \in [\knot_1, \knot_2); \\
0, & \text{otherwise}.
\end{cases}
\end{equation}

\paragraph{If \( p > 0 \), then}
\begin{equation}
B[\knotvector](x) = \frac{x - \knot_1}{\knot_{p+1} - \knot_1}
B[\knotvector^{-}](x) + \frac{\knot_{p+2} - x}{\knot_{p+2} - \knot_{2}}
B[\knotvector^{+}](x),
\end{equation}
where \( \knotvector^{+} \) and \( \knotvector^{-} \) are obtained by
dropping the first and last elements of \( \knotvector \) respectively:
\begin{align}
\knotvector^{+} = (\knot_2, \ldots, \knot_{p+2}), \quad \knotvector^{-} = (\knot_1, \ldots, \knot_{p+1}).
\end{align} 
In the cases of a vanishing denominator, the whole term is taken to be zero.
\end{definition}

Such univariate splines can be easily extended to higher dimensions through a
tensor product construction.

\begin{definition}
Let the polynomial multi-degree \( p = (p_1, \ldots, p_d ) \) and the \( d \)
local knot vectors \(\knotvector^{1}, \ldots, \knotvector^d \) be given. The
\define{\(d\)-variate tensor product B-spline} \( 
\bspline \colon \R^d \to \R \) is then defined as
\begin{equation}
\bspline(\vec{x}) = \prod_{i=1}^d B[\knotvector^i](x_i),
\end{equation}
where \( \vec{x} = (x_1, \ldots, x_d) \).
\end{definition}
The \define{support} of \( \bspline \) is the closure of the area where the
B-spline takes non-zero values, which we denote by:
\begin{equation}
\Supp(\bspline) = \closure{\set{\vec{x} \in \R^d : \bspline(\vec{x}) \neq 0}}.
\end{equation}

Since our B-spline construction is inherently local, we need to know when a
tensor product B-spline has \emph{minimal support} with respect to some box
mesh.

\begin{definition}[Informal]
A B-spline \( B = \bspline \) has \define{support on \( (\mesh, \mu) \)} if
all the knot lines of \( B \) occurs as meshlines in \( \mesh \). We say
that \( B \) has \define{minimal support on \( (\mesh, \mu) \)} if in
addition, all the knot lines of \( B \) occur consecutively in \( (\mesh,
\mu) \).
\end{definition}

One of the central concepts we will be addressing in this paper is the
\emph{overloading} of elements. We make this precise in the following
definition.
\begin{definition}\label{def:overloading}
Let a box partition \( \boxpartition \) of a domain \( \Omega \) and a
polynomial multi-degree \(p = (p_1, \ldots, p_d ) \) be given. Assume that we
construct a set \( \bsplines \) of B-splines degree \( p \) over the mesh \(
\mesh \) corresponding to \( \boxpartition\). We say that an element \(
\element \) is \define{overloaded with respect to \( \bsplines \)} if the
number of B-splines with support on \( \element \) is larger than the
dimension of the corresponding space of polynomials over this element, namely
\begin{equation}
\dim(\polynomials_p(\element)) = \prod_{i=1}^d (p_i+1).
\end{equation}
\end{definition}

We now proceed to review the definitions of LR B-splines and THB-splines.

\subsection{Locally refined spline spaces}

In preparation for the following discussion, we will adopt the notational
convention as in \cite{johannessenSimilaritiesDifferencesClassical2015a} in
order to differentiate between the distinct types of basis functions. Depending
on the underlying box partition, some of these types may coincide.

\begin{center}
\begin{tabular}{lrr}
\toprule
Type & Basis & Function \\
\midrule
Tensor Product B-spline & \( \bsplines \) & \( B \) \\
Truncated Hierarchical B-spline & \( \thbsplines \) & \( H \) \\
LR B-spline & \( \lrsplines \) & \( L \) \\
\bottomrule
\end{tabular}
\end{center}

Furthermore, in this and the following sections we will be dealing with box
partitions and spline spaces in \( \R^2 \), unless otherwise stated.

\subsubsection{LR-splines}

Locally Refinable Splines (LRB or LR-splines) was introduced by
\cite{dokkenPolynomialSplinesLocally2013}. The LR-spline framework permits
the insertion of local splits in a tensor product mesh, and subsequently
enables local refinement of the mesh. Being scaled tensor product B-splines,
LR-splines admit a set of nice properties. The set of LR B-splines form a
partition of unity. Their scaling weights are positive, meaning that they
satisfy the convex hull property, and are therefore inherently stable in
computations. Moreover, with some restrictions on the refinement process,
linear independence of the resulting set of functions can be guaranteed. 

\begin{figure}[htpb]
\centering
\subfloat[]{\includegraphics[width=0.49\linewidth]{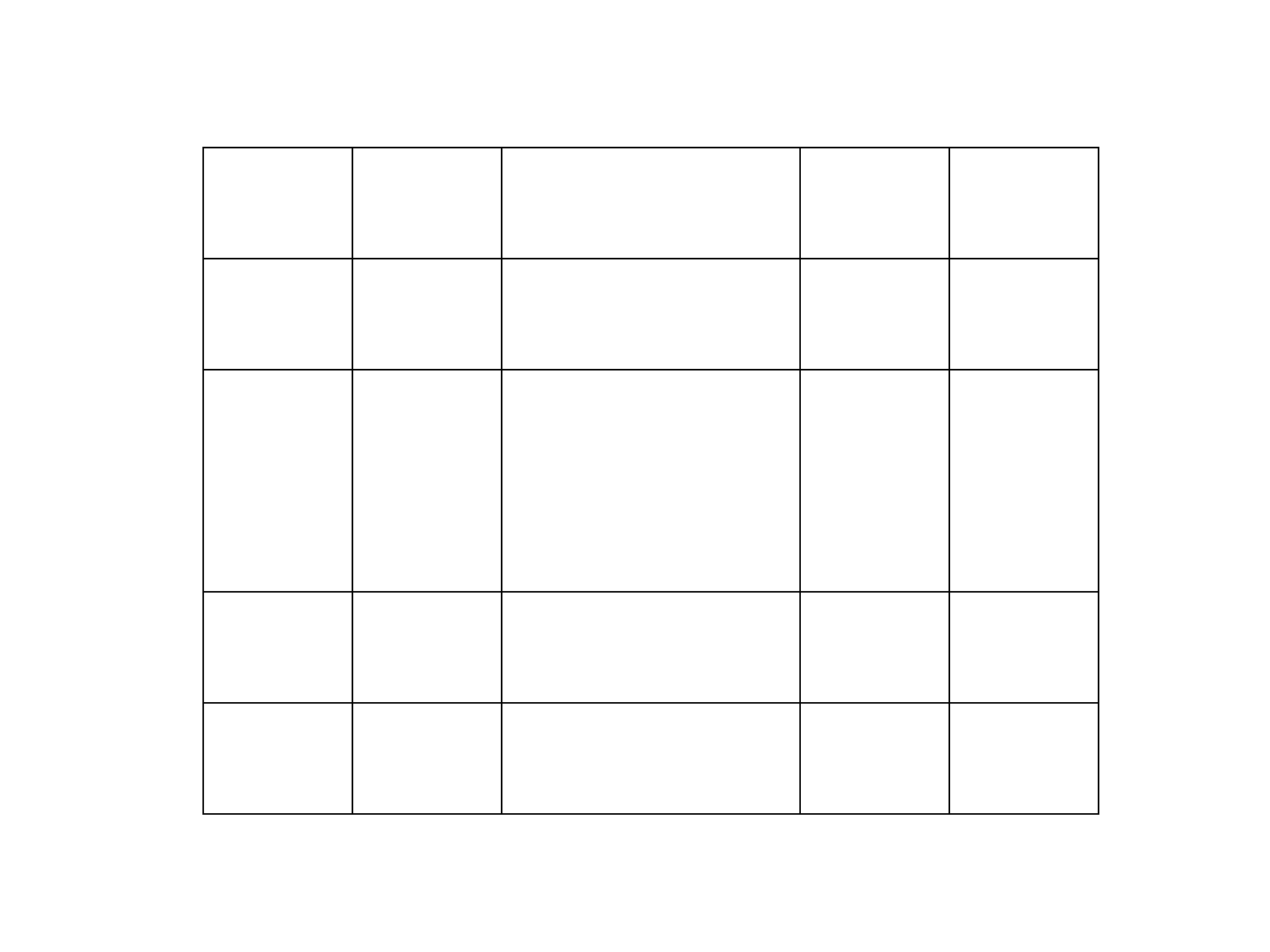}\label{sub:tensor_mesh}}
\subfloat[]{\includegraphics[width=0.49\linewidth]{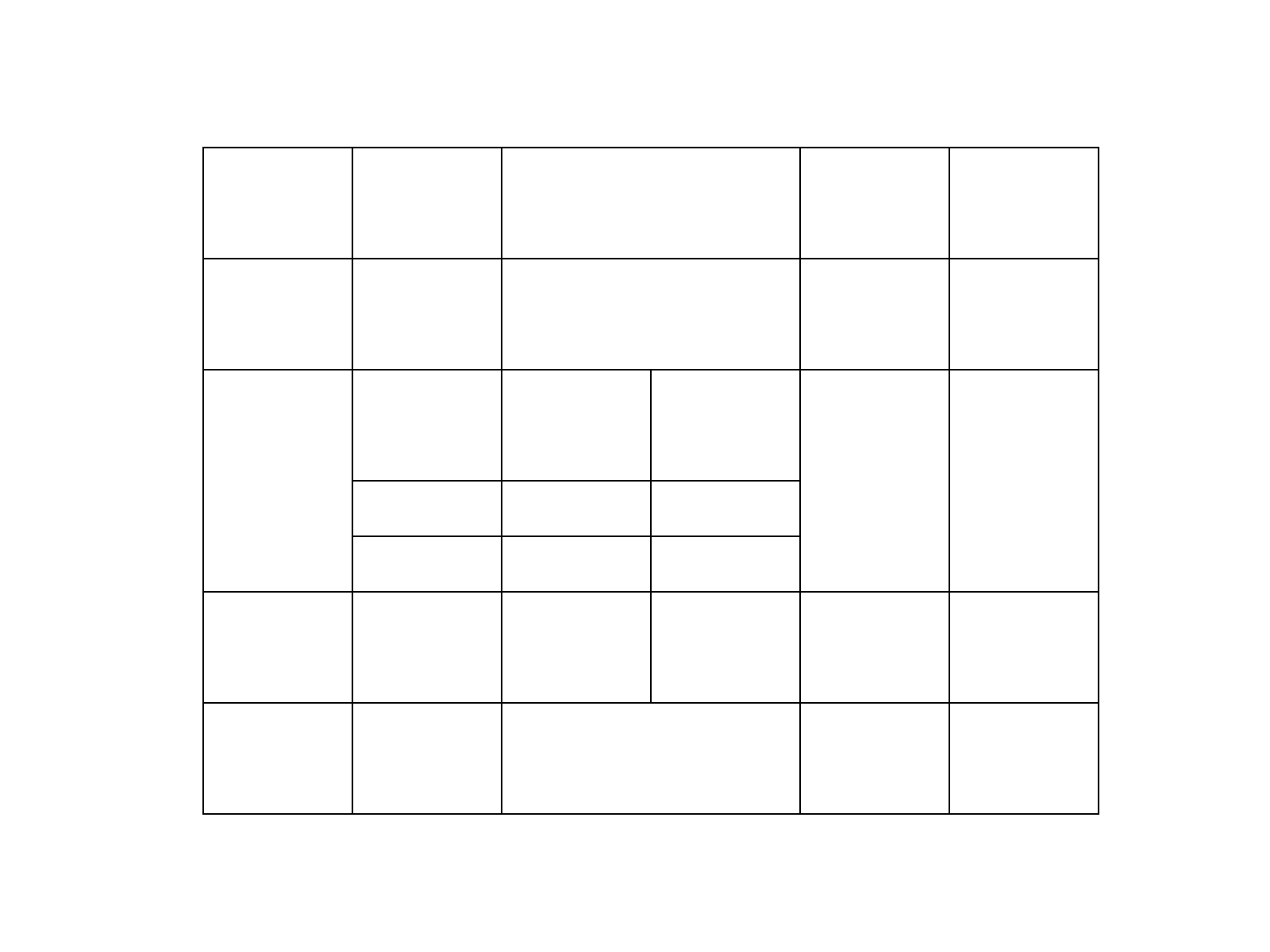}\label{sub:lr_mesh}}
\caption[An example of an LR-mesh]{In \protect\subref{sub:tensor_mesh}, an initial tensor product
mesh, which is also an LR-mesh. In \protect\subref{sub:lr_mesh}, an LR-mesh
obtained from the insertion of three meshlines in the initial tensor product
mesh from \protect\subref{sub:tensor_mesh}.}
\label{fig:lr_mesh}
\end{figure}

LR-splines are defined over so-called \emph{LR-meshes}, being special box
partitions. Starting from an initial tensor product mesh, meshlines are
inserted sequentially, yielding a sequence of box-meshes, where no meshline is
allowed to terminate in the middle of an element. This is formalized in the
following definition, and \cref{fig:lr_mesh} gives an example.

\begin{definition}
An \define{LR mesh} is a box mesh \( \mesh = \mesh_N\) resulting from a
sequence of meshline insertions in an initial tensor product mesh \( \mesh_1
\). That is 
\begin{equation}
\mesh_{i+1}= \mesh_{i} + \meshline_i
\end{equation}
for \( i = 1, \ldots, N-1 \) where each intermediate mesh is a box mesh. 
\end{definition}
\begin{remark}
We often think of an LR-mesh as a \emph{sequence} of intermediate meshes
\begin{equation}
\mesh = \mesh_N \supseteq \mesh_{N-1} \supseteq \cdots \supseteq \mesh_{2} \subseteq \mesh_{1}
\end{equation}
as each intermediate step is needed for the LR B-spline construction.
\end{remark}

Over such an LR-mesh we may define the associated set of LR B-splines algorithmically. Starting from an initial space of tensor product B-splines,
meshlines are inserted sequentially. Whenever a meshline completely traverses
the support of a B-spline, the B-spline is \emph{split} according to the knot
insertion procedure, and two new B-splines are added. The B-spline that was
split is removed.
\begin{definition}
Let \( \mesh \) be an LR-mesh over a domain \( \Omega \) and \( p = (p_1,
p_2) \) a polynomial bi-degree. We define the set \( \lrsplines(\mesh) \) of
\define{LR B-splines of degree \(p\) over \( \mesh \)} algorithmically as in
\vref{alg:lr_definition}.
\begin{algorithm}
\caption{The LR B-spline construction.}\label{alg:lr_definition}
\begin{algorithmic}
\State Let \( \lrsplines_1 \coloneqq \bsplines(\mesh_1) \) be the set of
tensor product B-splines on \( \mesh_1 \).
\For{each intermediate mesh \(\mesh_{i+1} = \mesh_{i} + \meshline_i \), with \( i = 1, \ldots, N - 1 \)}
\State{\(\lrsplines_{i+1} \coloneqq \lrsplines_{i}\)}
\While{there exists \( B \in \lrsplines_{i+1}\) without \emph{minimal
support} on \( \mesh_{i+1} \)}
\State \( B^+, B^- = \Call{Split}{B, \meshline_i} \) \Comment{knotline
insertion}
\State \( \lrsplines_{i+1} = (\lrsplines_{i+1} \setminus \set{B}) \cup \set{B^{+}, B^{-}} \) \Comment{update the set of B-splines}
\EndWhile
\EndFor
\State{\(\lrsplines(\mesh) \coloneqq \lrsplines_N\)}
\end{algorithmic} 
\end{algorithm}
\end{definition}

\begin{remark}
Note that all LR B-splines have \emph{minimal support} on the resulting mesh.
This is by construction. However, there is an important distinction to be made,
namely that the set of LR B-splines \emph{differ} from the set of minimal
support B-splines on the resulting mesh. This is due to the LR refinement
procedure putting some restrictions on the resulting mesh. A survey on the
properties of LR-splines and minimal support B-splines are given in
\cite{patriziLinearDependenceBivariate2018}.
\end{remark}

\subsubsection{Truncated hierarchical B-splines}

Hierarchical B-splines, first introduced in
\cite{forseyHierarchicalBsplineRefinement1988}, is a method for specifying
locally refined spline spaces on hierarchical meshes. Recall that a
hierarchical mesh consists of regions corresponding to various levels of tensor
product grids. The hierarchical B-spline construction involves replacing any
B-spline with support completely contained in a region of a finer level by
B-splines at this finer level. This procedure will, however, lead to coarse
B-splines partially overlapping the finer regions, and does not constitute a
partition of unity. 

A remedy to this problem came with the introduction of Truncated Hierarchical
B-splines \cite{THBsplinesTruncatedBasis2012}, where B-splines on a coarse
level is \emph{truncated} by B-splines on a finer level. This leads to the
resulting set of B-splines forming a partition of unity. The construction
relies on the \emph{truncation operator}. Recall that a spline \( f \in
\SPAN(\bsplines_\ell) \) can be represented in terms of the finer basis \(
\bsplines_{\ell + 1} \):
\begin{equation}
f = \sum_{\mathclap{B_i \in \bsplines_{\ell + 1}}} c^{\ell + 1}_i(f) B_i,
\end{equation}
where \( c^{\ell + 1}_i(f) \) is the coefficient multiplying \( B_i \) in the
representation of \( f \) in terms of \( \bsplines_{\ell + 1}\).
For uniform B-splines, this relation is often called the \emph{two-scale
relation}. The truncation operator is defined as follows:

\begin{definition}
Let \( B \in \bsplines_\ell \) be a coarse B-spline. The
\define{truncation} with respect to the set of fine B-splines \(
\bsplines_{\ell + 1} \) and the corresponding region \( \Omega_{\ell + 1}
\) is
\begin{equation}
\trunc^{\ell + 1} B \coloneqq \sum_{\mathclap{\substack{B_i \in
\bsplines_{\ell+1} \\ \Supp B_i \not\subseteq \Omega_{\ell + 1}}}}
c^{\ell + 1}_i(B)B_i.
\end{equation}
\end{definition}

\begin{remark}
The definition above represents the truncation operator in an
\emph{additive} sense, where the contributions from the finer level are
summed up. It is also possible to represent the truncation operator
\emph{subtractively}, by instead removing the bits of the representation
that have been replaced by finer B-splines: 
\begin{equation}
\trunc^{\ell + 1} B = B - \sum_{\mathclap{\substack{B_i \in
\bsplines_{\ell+1} \\ \Supp{B_i} \subseteq \Omega_{\ell + 1}}}} c^{\ell
+ 1}_i(B) B_i
\end{equation}
\end{remark}

\begin{definition}
Let \( \mesh \) be a hierarchical mesh over a domain \( \Omega \) (see
\cref{def:hierarchical_mesh}) and \( p = (p_1, p_2) \) a polynomial bi-degree.
On each level \( \ell = 1, \ldots, N \) we have a tensor product spline space
\( V_\ell \) spanned by a collection of B-splines \( \bsplines_\ell =
\bsplines(\mesh_\ell)\). We define the set of \define{THB-splines of degree \(
p \) over \( \mesh \)} algorithmically as in \cref{alg:thb_definition}.

\begin{algorithm}
\caption{The THB-spline construction.}\label{alg:thb_definition}
\begin{algorithmic}
\State Let \( \thbsplines_1 = \bsplines(\mesh_1) \) be the set of tensor product B-splines on \( \mesh_1 \).
\For{each level \( \ell = 1, \ldots, N-1 \)}
\State{
\begin{flalign*}
\thbsplines_{\ell + 1}^{\mathrm{trunc}} &\coloneqq \set{\trunc^{\ell + 1} H : H \in \thbsplines_{\ell} \text{ and } \Supp(H) \not\subseteq \Omega_{\ell+1}} \\
\thbsplines_{\ell + 1}^{\mathrm{new}} &\coloneqq 
\set{B \in \bsplines_{\ell+1} : \Supp(\beta) \subseteq \Omega_{\ell + 1}} \\
\thbsplines_{\ell + 1} &\coloneqq \thbsplines_{\ell+1}^{\mathrm{trunc}} \cup \thbsplines_{\ell + 1}^{\mathrm{new}} 
\end{flalign*}}
\EndFor
\State \( \thbsplines = \thbsplines_{N} \).
\end{algorithmic} 
\end{algorithm}
\end{definition}

\begin{remark}
Note that in the cases where a B-spline \( B \in \bsplines_{\ell} \) is
truncated with respect to \( \bsplines_{\ell+1} \) and \( \Omega_{\ell + 1}
\) and the support of \( B \) happen to be entirely contained in \(
\Omega_{\ell + 1} \), the truncation operator completely removes the coarse
B-spline. In the THB-spline construction, this has the effect of replacing
the coarse B-splines with fine B-splines defined in its support.
\end{remark}

\begin{remark}
A simple framework for the implementation of truncated hierarchical
B-splines is given in \cite{garauAlgorithmsImplementationAdaptive2018a},
and this serves as a good introduction to the many ways such splines have
been implemented in the literature. Efficient algorithms for the assembly
of finite element matrices are also presented.
\end{remark}

\subsubsection{T-splines}

\begin{figure}[htpb]
\centering
\subfloat[]{\includegraphics[width=0.49\linewidth, height=0.49\linewidth]{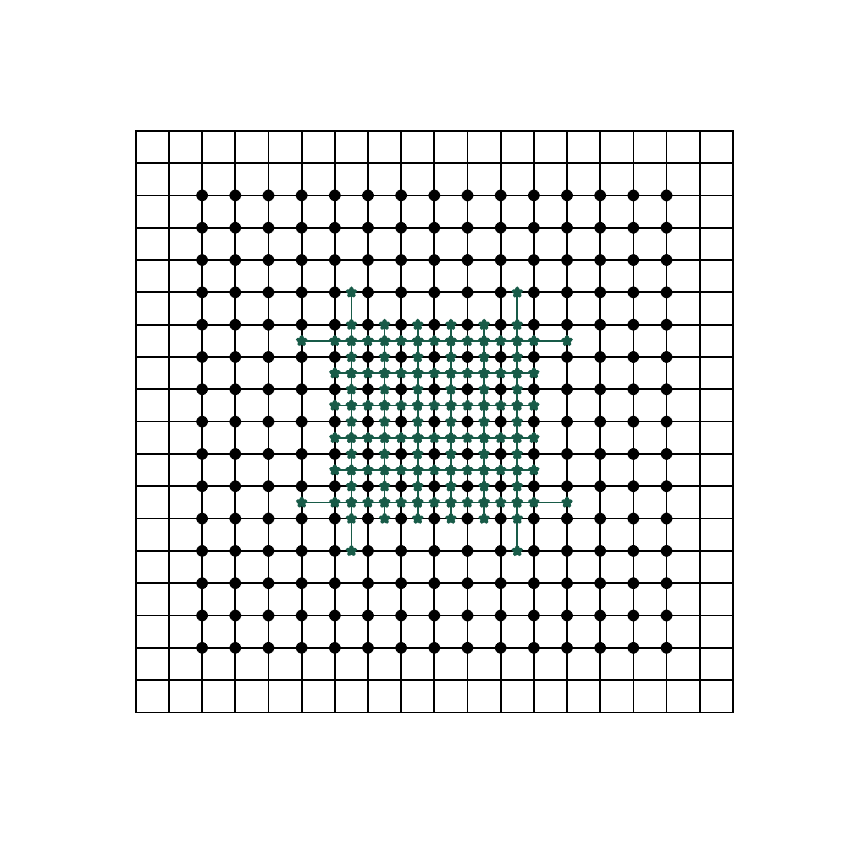}\label{sub:t_mesh}}
\subfloat[]{\includegraphics[width=0.49\linewidth, height=0.49\linewidth]{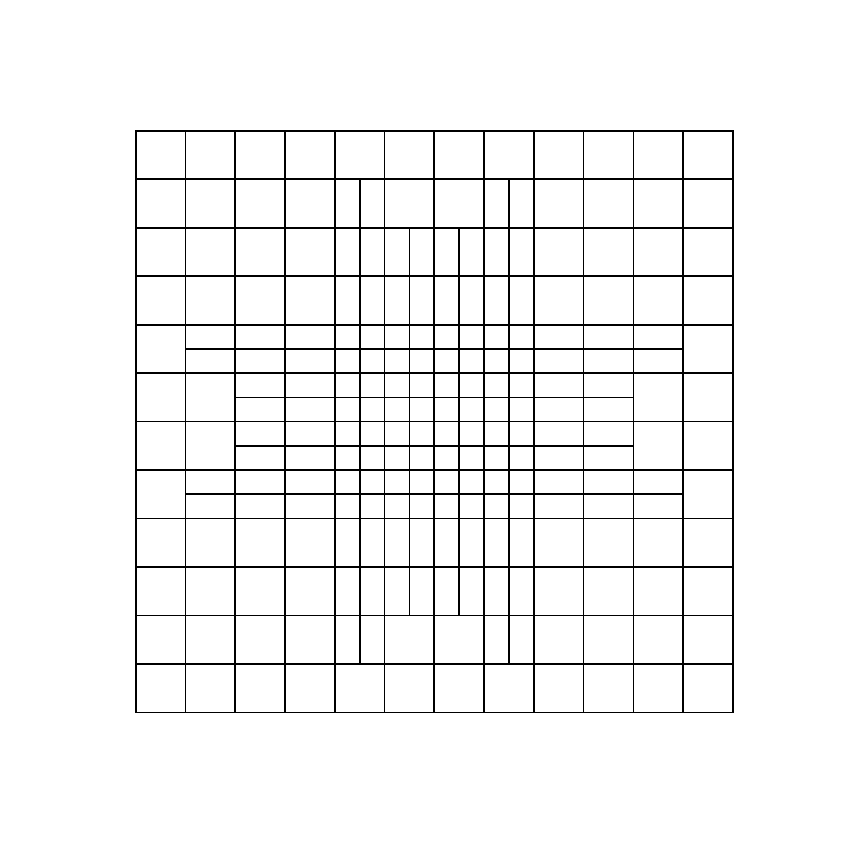}\label{sub:t_mesh_lr}}
\caption[A T-spline mesh]{In \protect\subref{sub:t_mesh} a T-spline mesh in
the index domain. The dots denote 
Greville points or ``anchors'' for each
individual B-spline. A black dot is a B-spline at level \( \ell = 0 \) and
a green star a B-spline at level \( \ell = 1\). The resulting spline space can be replicated
by an LR-mesh without overloaded elements (c.f.
\cref{sub:overloading_lrb-tno}), as displayed in
\protect\subref{sub:t_mesh_lr}. Here we have used multiplicity \(m = 4\) along
the boundary.} 
\label{fig:t_mesh}
\end{figure}

While not directly addressed in this paper, we briefly mention T-splines as
LRB with local modifications to the LR-meshes used in this paper is able
to reproduce the spline space generated by semi-standard T-splines \cite{Sederberg-2004}
and Analysis Suitable T-splines \cite{Scott-2012}. An example of an Analysis
Suitable T-mesh in the index domain is displayed in \cref{fig:t_mesh} to the
left, with the corresponding LR-mesh to the right. This is a close up of the
structure of a mesh similar to the one used in Figure
\ref{sub:overloading_lrb-tno}.

%% file: local_modifications_and_overloading.tex
\section{Local modification of meshes and the reduction of overloading}
\label{sec:local_modifications_and_overloading}

In this section we take a deeper look at the overloading of elements, and how
local modifications to the mesh may be used to remedy this. Recall from the
previous definition that an element \( \element \) in a box-partition \(
\boxpartition \) is said to be \emph{overloaded} if the number of supported
B-splines on the element exceeds the number needed to span the full polynomial
space over the element. 

\begin{figure}[htpb]
	\centering
	\subfloat[LRB]{\includegraphics[width=1\linewidth]{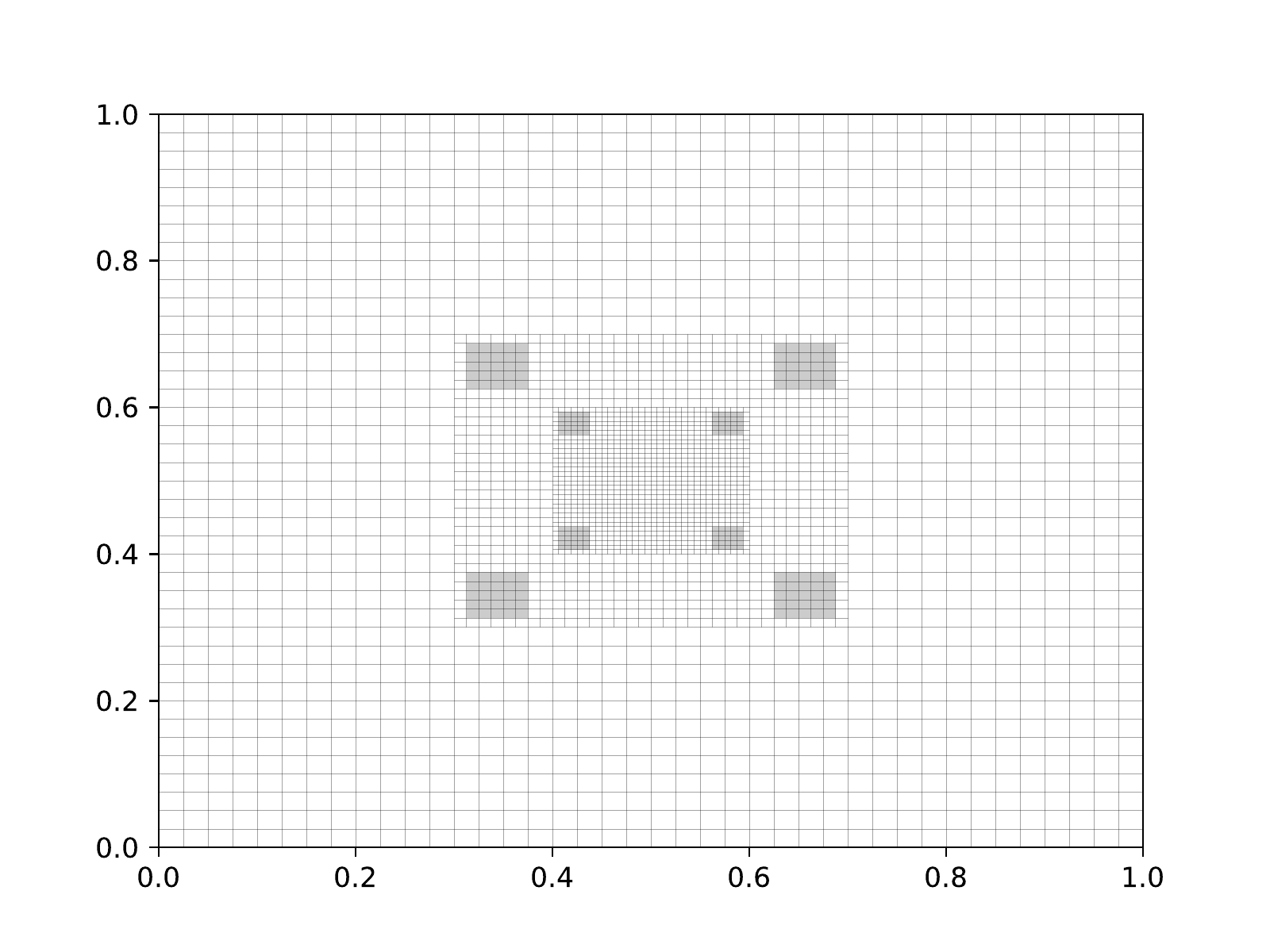}\label{sub:overloading_lrb}}

	\subfloat[THB]{\includegraphics[width=1\linewidth]{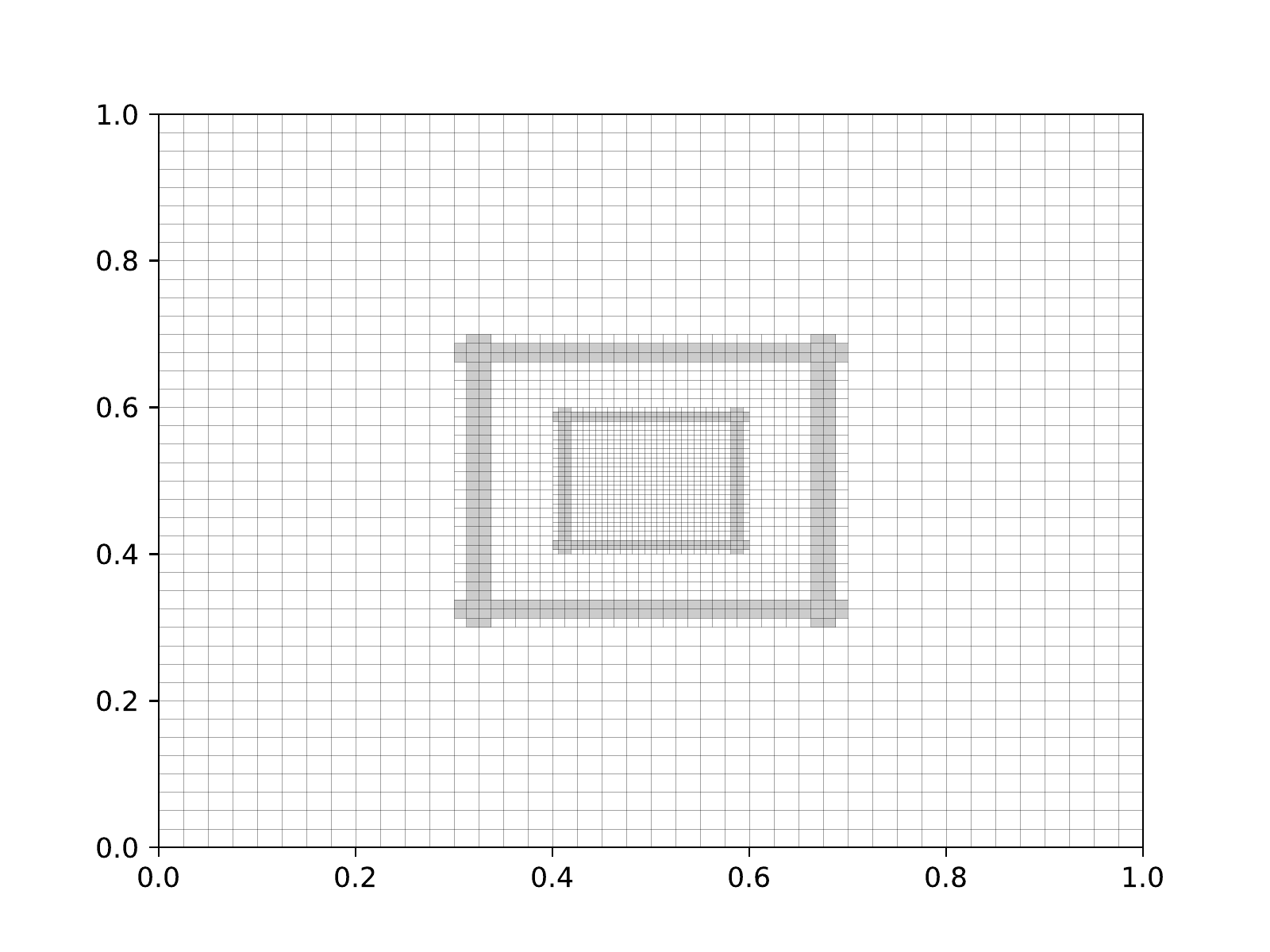}\label{sub:overloading_thb}}
  \caption[Overloading patterns for THB and unmodified LRB]{The overloading patterns on a hierarchical mesh with two
		levels of refinement. In \protect\subref{sub:overloading_lrb},
		we see that regions in the corners of the refined regions are
		overloaded, due to the influence of four LR B-splines from
		the coarser layer, whose support has not been split by any
    newly introduced meshlines. In \protect\subref{sub:overloading_thb} we
    observe ``bands'' of overloaded elements along the boundary between two
  consecutive refinement levels for THB, arising due to the fact that fine
  B-splines must be completely contained in the support of a coarse B-spline
  before truncation occurs.}
	\label{fig:overloading_example}
\end{figure}

We are interested in such overloaded regions, because by reducing or removing
completely the overloading on elements we may
\begin{enumerate}
	\item reduce the bandwidth of the resulting finite element matrices; and 
	\item improve conditioning of finite element matrices.
\end{enumerate}
Such overloaded regions occur for LRB in convex corners of a fine hierarchical
level, where a large B-spline from one hierarchical level overlaps several
elements of a finer hierarchical level. For THB, overloading occurs along any
border between two hierarchical levels. By coloring in elements with too many
supported B-splines we obtain a visualization of this phenomenon, as seen in
\cref{fig:overloading_example} on a hierarchical mesh with three levels of
refinement.

\begin{figure}[htpb]
	\centering
	\subfloat[LRBNO]{\includegraphics[width=1\linewidth]{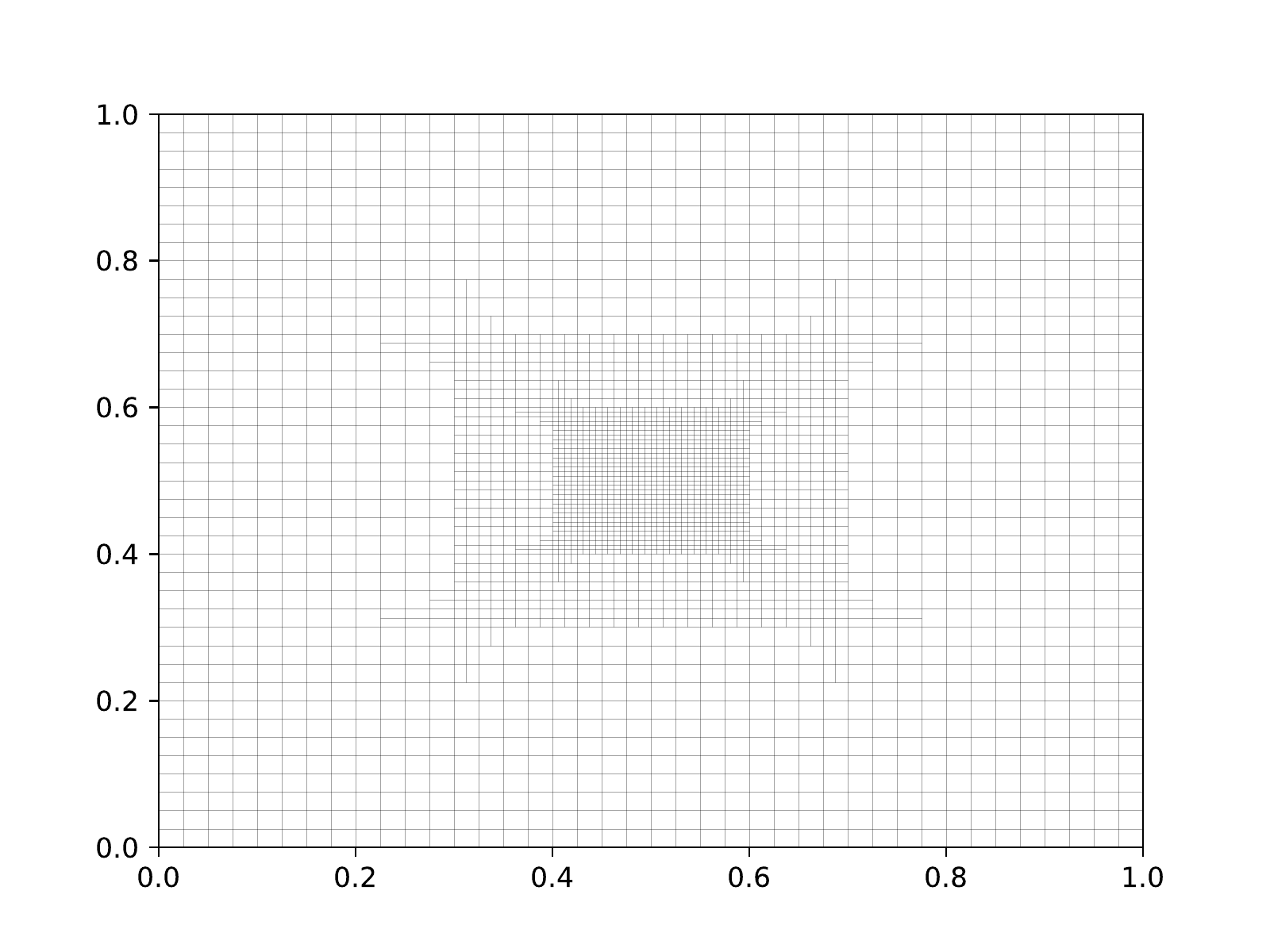}\label{sub:overloading_lrb-no}}

	\subfloat[T-LRBNO]{\includegraphics[width=1\linewidth]{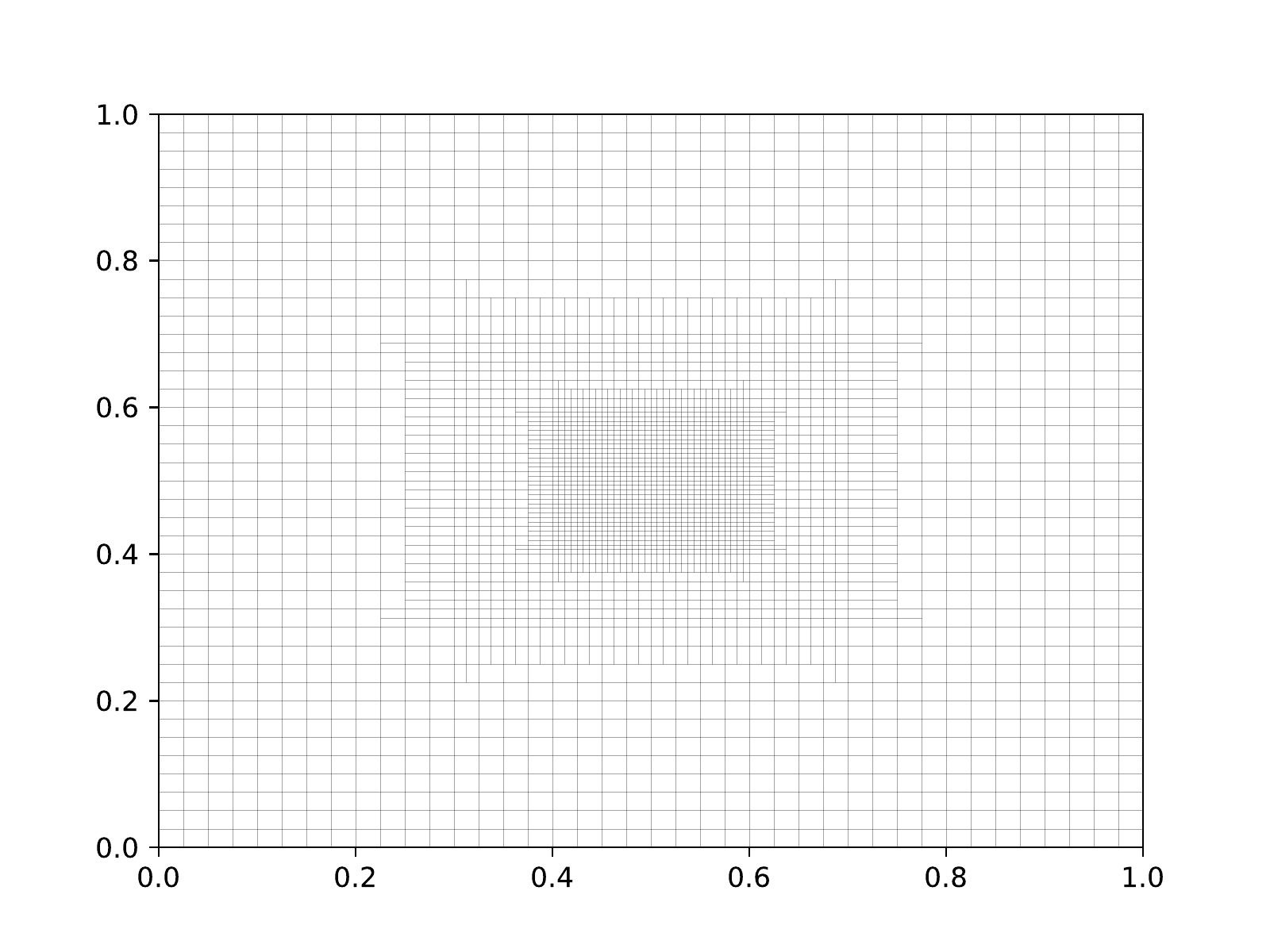}\label{sub:overloading_lrb-tno}}
  \caption[Overloading patterns for LRB with local modifications]{Two different local modifications with the effect of
		completely removing the overloaded elements. In
		\protect\subref{sub:overloading_lrb-no} we extend the meshlines
		closest to the convex corner by three elements, and the
		meshline next to them by one element. This has the effect of
		completely removing the overloading on the corner elements. 
    In \protect\subref{sub:overloading_lrb-tno}, we make a mesh that can be
  defined using T-splines that has no overloading. As in (a) meshlines closest
to the corners are extended by three, while meshlines at the borders between
refinement levels are extended by two as in Figure \ref{fig:t_mesh}.}
	\label{fig:overloading_example_modified}
\end{figure}

In order to reduce, or completely remove such overloaded regions, we may for
LRB extend meshlines from the fine hierarchical level to the coarse level, in
order to split the culprit B-splines. The length needed for this extended
meshline depends on the polynomial degree of the B-spline to be split. In
\cref{fig:overloading_example_modified} we see the effects of three types of
meshline extension to the LRB-mesh from \cref{fig:overloading_example} for a
space of bi-cubic splines.

\begin{figure}
  \centering
  \subfloat{\includegraphics[width=0.49\linewidth]{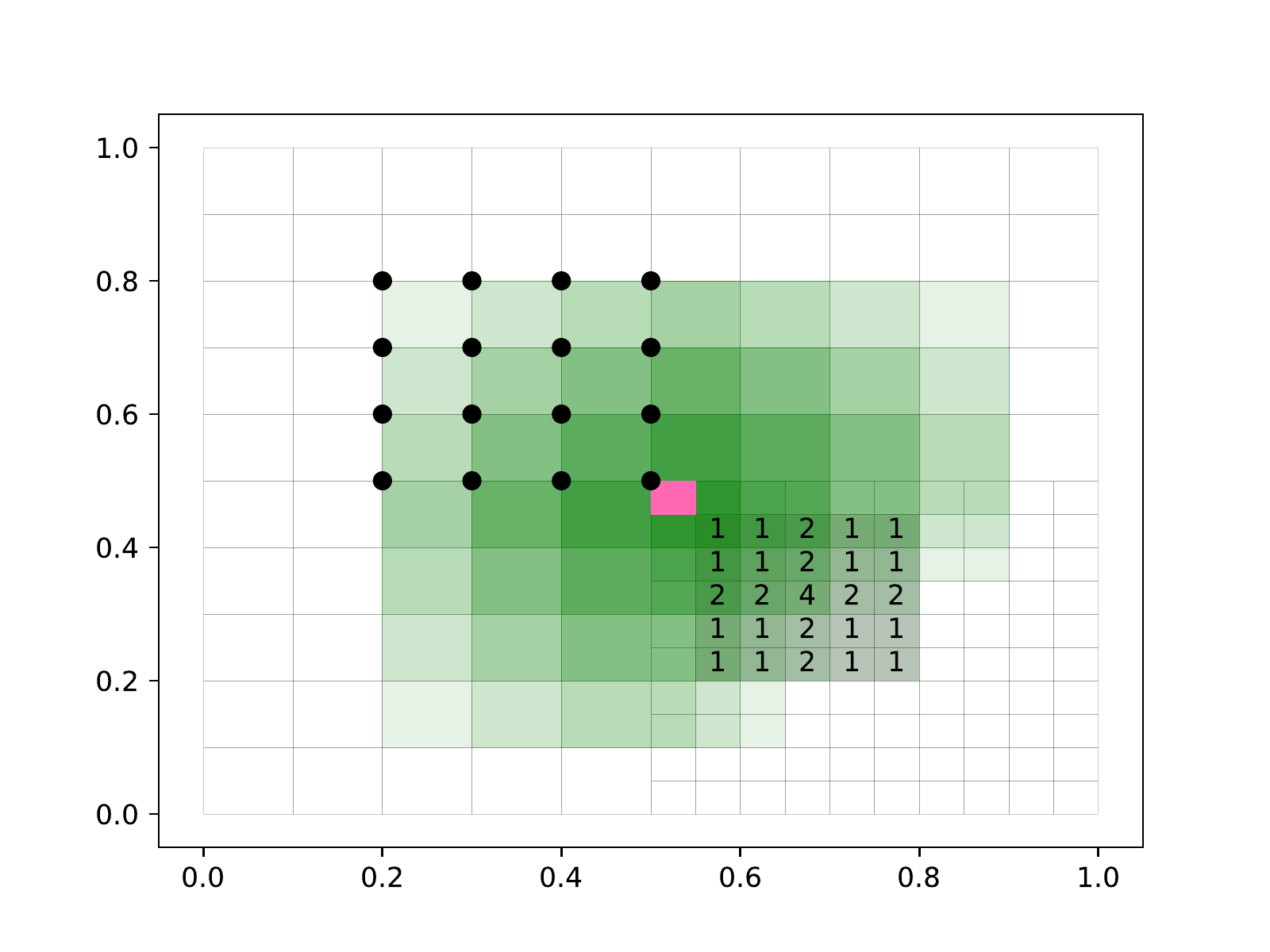}}

  \subfloat{\includegraphics[width=0.49\linewidth]{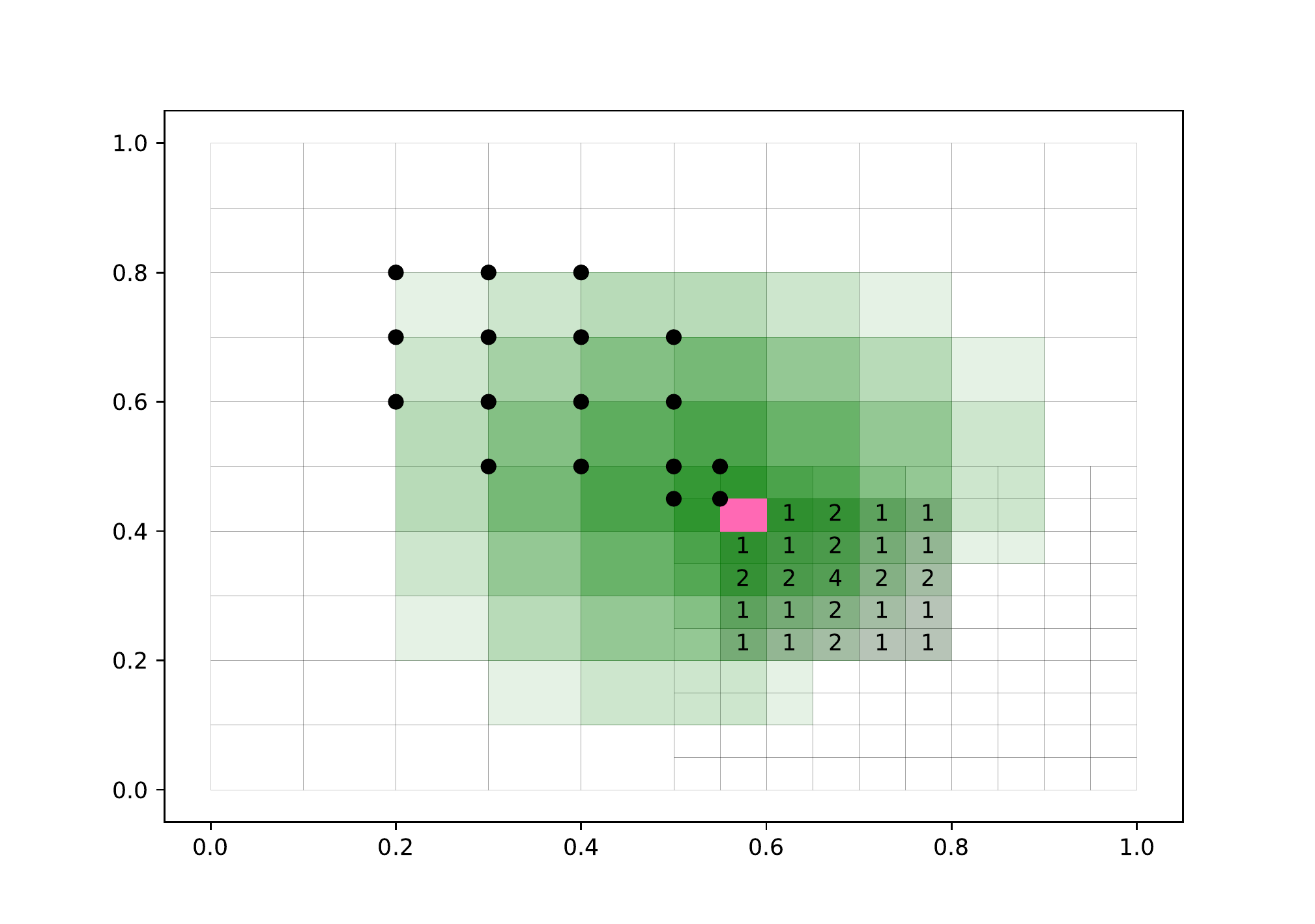}}
  \subfloat{\includegraphics[width=0.49\linewidth]{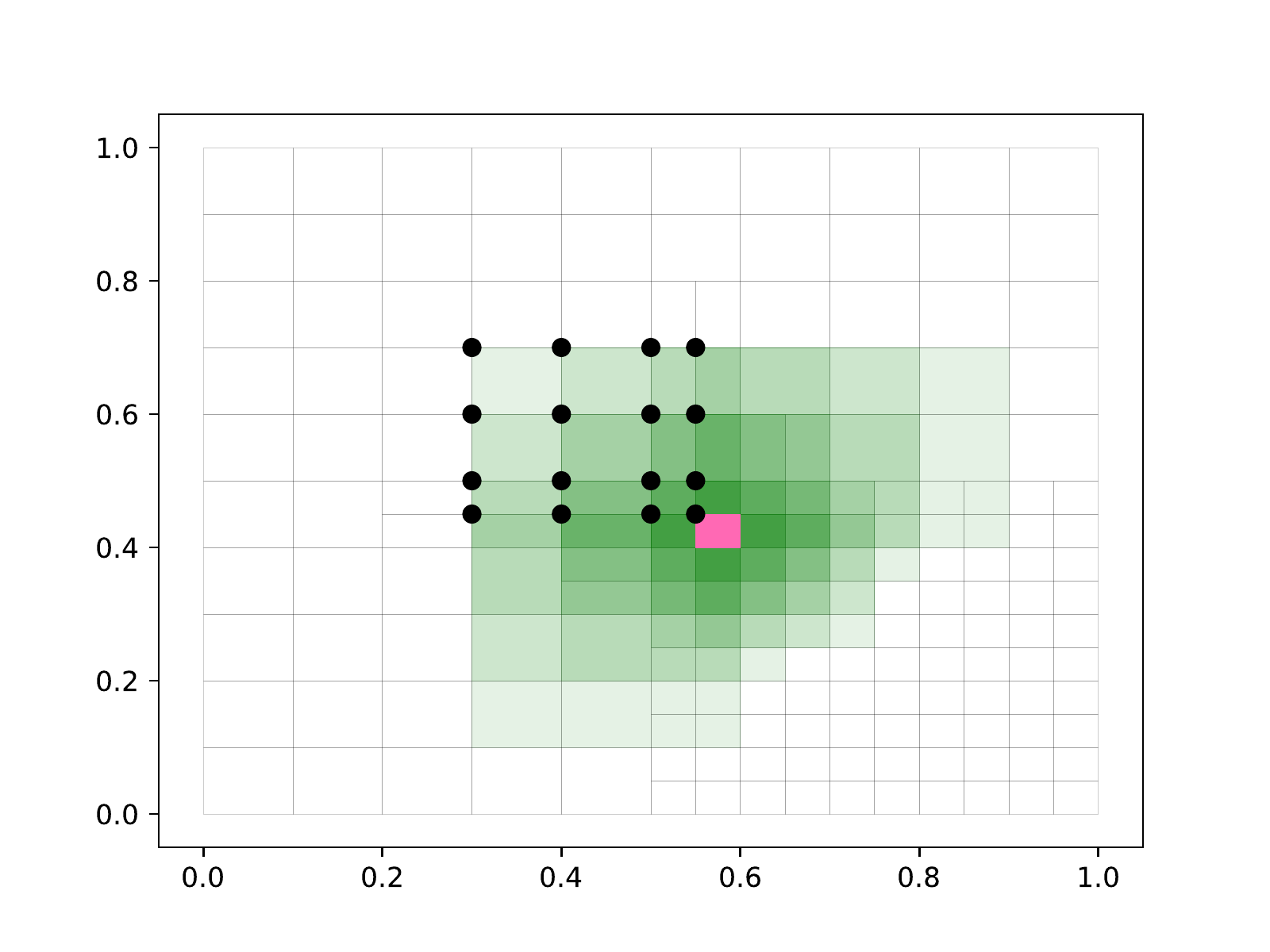}}

  \subfloat{\includegraphics[width=0.49\linewidth]{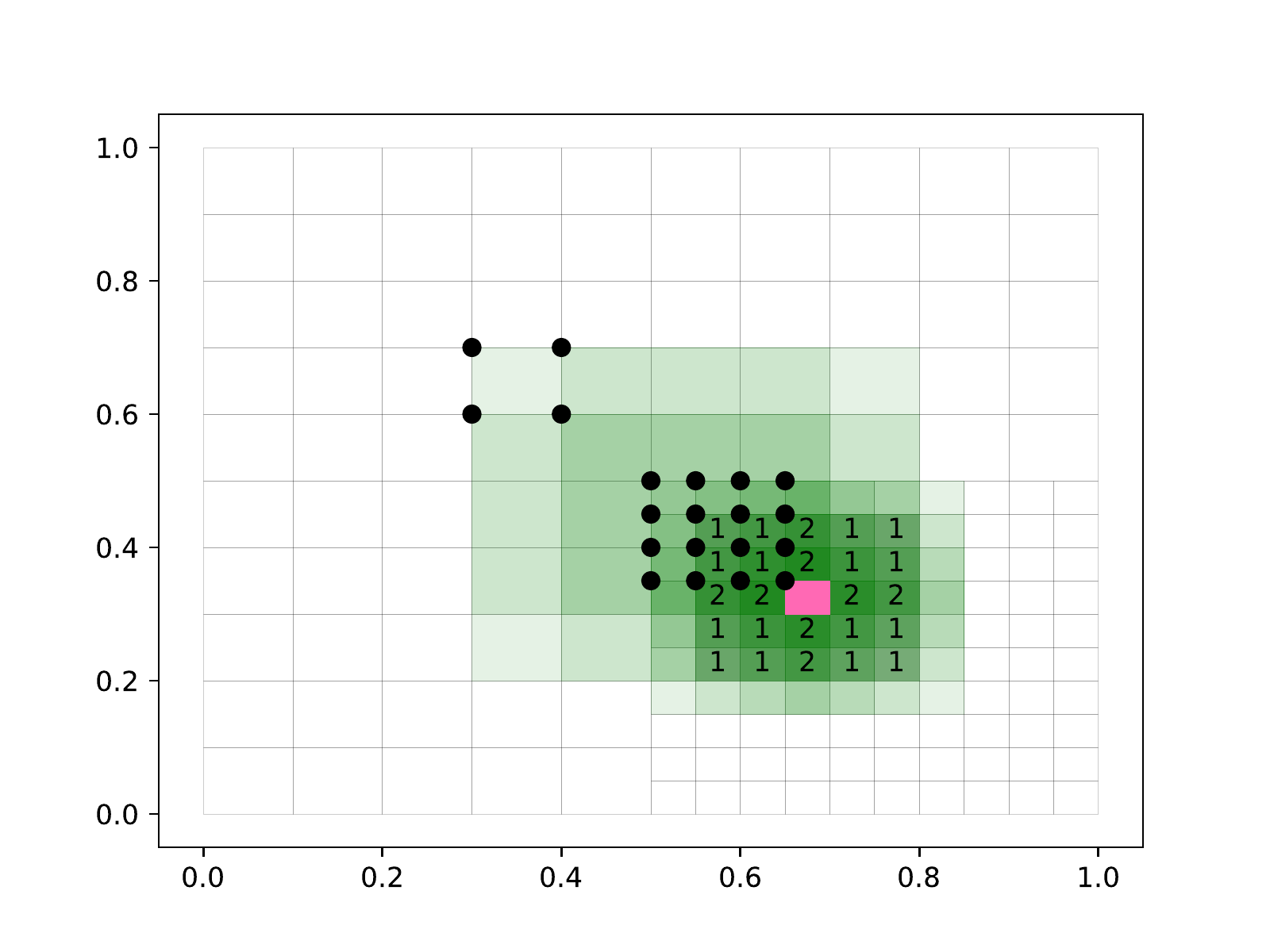}}
  \subfloat{\includegraphics[width=0.49\linewidth]{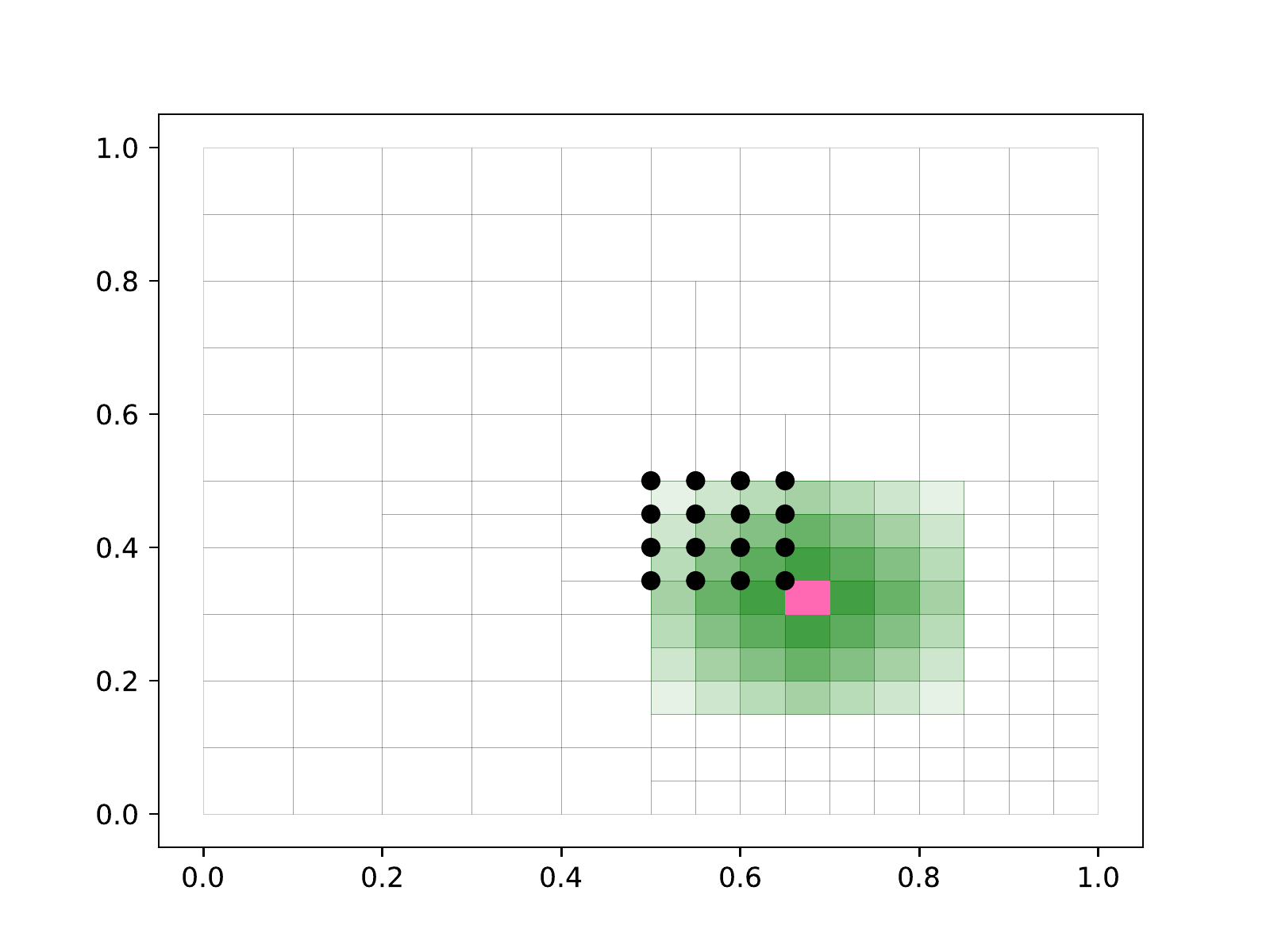}}
  \caption{
    The effects of extending meshlines on the bi-cubic B-splines
    covering the element in pink.  The upper left corner of each B-spline is
    marked with a black dot. The knotlines of each B-spline can be identified
    by starting from the dot and going four knotlines to the right/down. We
    chose to not use Greville points as some overloaded configurations produce
    overlapping Greville points.  In the \emph{upper mesh} we look at the
    element just inside the corner of the region refined, and no overloading
    occurs.  In the \emph{middle meshes} we move one element diagonally into
    the refined region. Before refinement the overloading is one, and after 
    additional lines are inserted the overloading is removed.  In the
    \emph{bottom meshes} we move two additional element diagonally into the
    refined region, Before refinement the overloading is four, after 
    additional lines are inserted the overloading is removed.
  }
  \label{fig:overload_reduction}
\end{figure}

In order to capture what is happening, we take a closer look at overloading in a
convex corner in \cref{fig:overload_reduction} where we show how B-splines from
the coarse level of a hierarchical mesh may overlap with B-splines from the
fine level in such a way that too many B-splines are active over a given
element.

%% file: numerical_experiments.tex
\section{Numerical experiments}
\label{sec:numerical_experiments}

In order to compare the methods addressed in this paper, we assemble the mass
and stiffness matrices associated to discretizations of partial differential
equations using IgA or FEM. By computing the condition number of these
matrices, we obtain a metric useful for comparison. These matrices arise
amongst others in the discretizations of the Poisson equation, and in the
computation of the \(L^2\)-projection of a function.

\subsection{\(L^2\)-projection}

Given a domain \( \Omega \), a function \( f \colon \Omega \to \R \) in some
space of functions \( V \), and a finite-dimensional subspace \( V_h \) of \( V
\), we are interested in finding the function \( u \in V_h \) that minimizes
the \( L^2(\Omega)\)-error
\begin{equation}
\| e\|_{L^2} = \| u - f\|_{L^2}.
\end{equation}
This can be reformulated as a variational equation by requiring \( s \) to
satisfy
\begin{equation}
\int_\Omega u v \dif{\Omega} = \int_\Omega f v \dif{\Omega},
\end{equation}
for all \( v \in V_h \). By introducing a basis \( \set{\varphi_1, \ldots,
\varphi_N} \) for \( V_h \), which in our case will be one of the THB or
LRB-bases, we may write this as a linear equation
\begin{equation}
\mat{M}\vec{c} = \vec{b},
\end{equation}
where \( \mat{M} \) is the \emph{mass matrix} and \( \vec{c} \) is the vector
of coefficients representing \( u \) in our chosen basis. The entries for \(
\mat{M} \) and the right-hand side \( \vec{b} \) are given as
\begin{equation}
\mat{M}_{ij} = \int_\Omega \varphi_i \varphi_j \dif{\Omega}\,, \quad \vec{b}_j = \int_\Omega f \varphi_j \dif{\Omega}.
\end{equation}

\subsection{The Poisson equation}

A commonly encountered differential equation is the \emph{Poisson equation}.
Given a function \( f \colon \Omega \to \R \), we wish to find a function \( u
\) in a space of admissible functions \( V \) such that 
\begin{equation}
\laplace{u} = f \text{ in } \Gamma,
\end{equation}
subject to the boundary conditions
\begin{equation}
u = 0 \text{ on } \Gamma_D\,, \quad \pd{u}{\vec{n}} = g \text{ on } \Gamma_N.
\end{equation}
Here \( \Gamma_D \) denotes the \emph{Dirichlet}-boundary and \( \Gamma_N \)
the \emph{Neumann}-boundary. We assume \( \boundary\Omega = \Gamma_D \cup
\Gamma_N \) and \( \Gamma_D \cap \Gamma_N = \emptyset \). Furthermore, \(
\vec{n} \) is the outward facing boundary normal to \( \Omega \) and \( g \) is
the prescribed flux along the boundary.

As for the \(L^2\)-projection, the Poisson equation can be reformulated as a
variational equation by multiplying the equation by a test-function, and
integrating over the domain. The requirements on the smoothness of the sought
solution \( u \) can be relaxed, by moving some derivatives onto the
test-functions.  Again, we seek the
solution \( u \) in a subspace \(V_h\) of \( V \) spanned by a set of basis
functions \( \set{\varphi_1, \ldots, \varphi_N} \). The variational form of the
Poisson-equation reads
\begin{equation}
\int_\Omega \grad{u} \grad{v} \dif{\Omega} = \int_\Omega f v \dif{\Omega} - \int_{\Gamma_N} gv \dif{S},
\end{equation}
for all \( v \in V \). Rewriting this in terms of the basis functions, we
obtain the linear equation
\begin{equation}
\mat{A} \vec{c} = \vec{b},
\end{equation}
where \( \mat{A} \) is the \emph{stiffness matrix} of the problem. The entries
of \( \mat{A} \) and \( \vec{b} \) are given as
\begin{equation}
\mat{A}_{ij} = \int_\Omega \grad{\varphi_i} \cdot \grad{\varphi_j}
\dif{\Omega}\,, \quad \vec{b}_j = \int_\Omega f \varphi_j \dif{\Omega} -
\int_{\Gamma_N} \varphi_j g \dif{S}.
\end{equation}

\begin{figure}[htpb]
\centering
\includegraphics[width=1\linewidth]{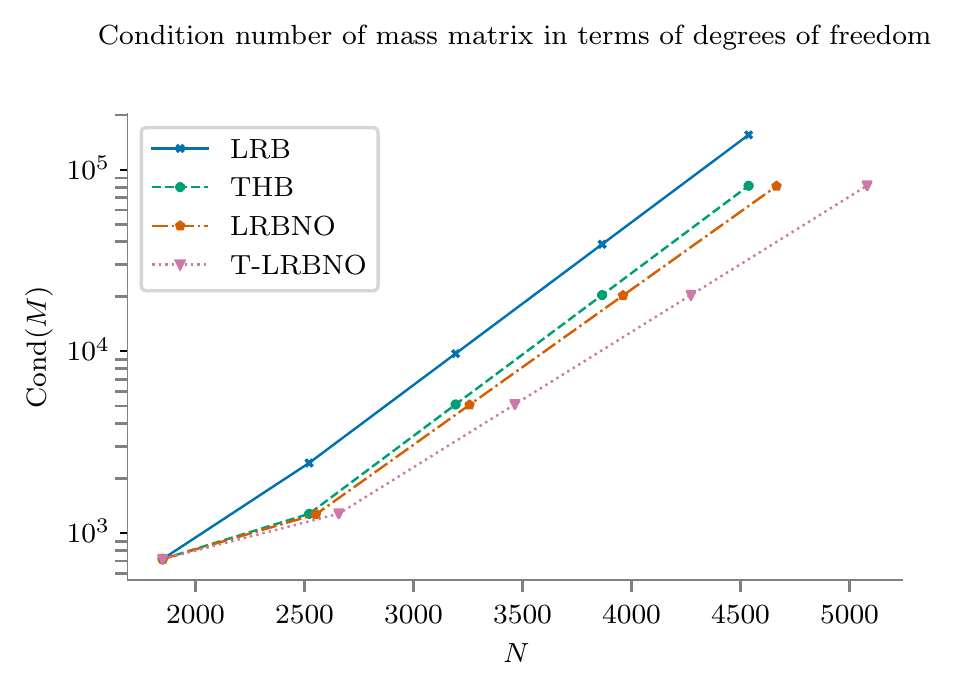}\label{sub:central_mass_ndof}
\caption{The condition numbers for mass matrices over a centrally refined
hierarchical mesh for six levels of refinement. In the figures \( N \)
denotes the number of degrees of freedom in the corresponding space.}
\label{fig:central_refinement_condition_mass}
\end{figure}

\begin{figure}[htpb]
\centering
\includegraphics[width=1\linewidth]{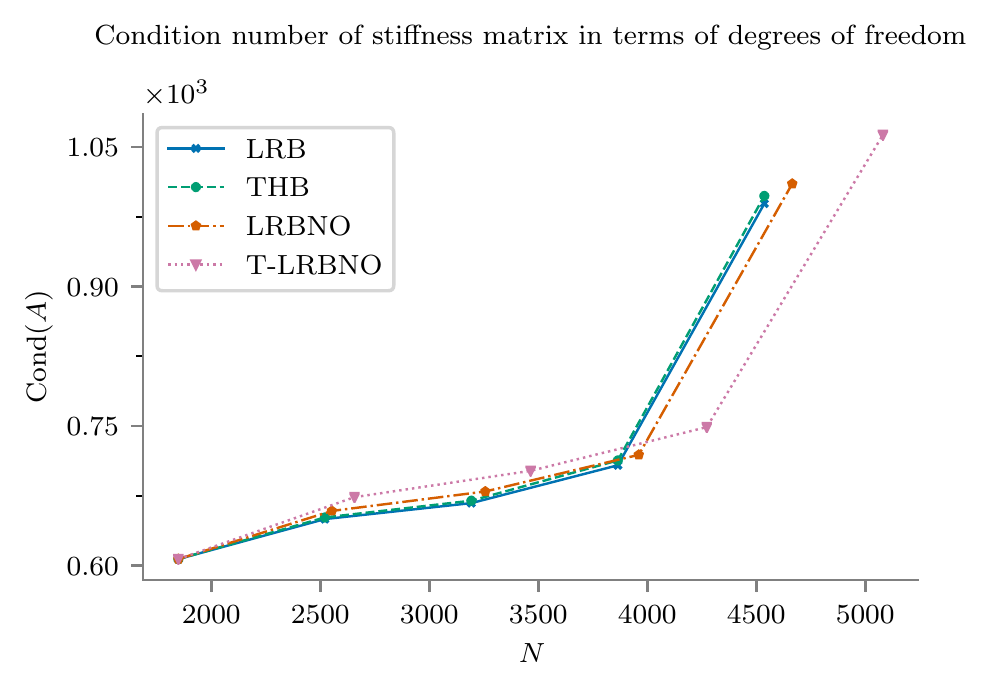}
\caption{The condition numbers for stiffness matrices
over a centrally refined hierarchical mesh for six levels of
refinement. In the figures \( N \) denotes the number of degrees of freedom
in the corresponding space.}
\label{fig:central_refinement_condition_stiffness}
\end{figure}

\subsection{Condition Numbers}

The \emph{condition number} of a matrix \( \mat{B} \in \R^{n\times n} \)
quantifies how sensitive the solution \( x \) to the linear system \(
\mat{B}\vec{x} = \vec{y} \) is to small perturbations both in \( \mat{B} \) and
the right-hand side \(\vec{y} \) and is formally defined as 
\begin{equation}
\cond{\mat{B}} \coloneqq \| \mat{B} \| \| \mat{B}^{-1} \|,
\end{equation}
where \( \| \cdot \| \) is some matrix norm. Note that the condition number is
norm-dependent, but all matrix norms are equivalent on \( \R^{n \times n} \).
We will be computing the condition numbers in the 2-norm, and in this specific
setting for \emph{normal} matrices the condition number can be computed as the
ratio between the largest and smallest eigenvalue
\begin{equation}
\cond{\mat{B}} = \frac{|\lambda_1(\mat{B})|}{|\lambda_n(\mat{B})|}.
\end{equation}
Here \( \lambda_1 \geq \lambda_2 \geq \ldots \geq \lambda_n \), i.e., ordered
in a decreasing fashion.

As in \cite{johannessenSimilaritiesDifferencesClassical2015a}, we chose to
estimate the condition numbers of the matrices \emph{before} imposing any
boundary conditions, as imposing boundary conditions can have a large impact on
the conditioning of the matrix. The mass matrix \( \mat{M} \) is non-singular,
even with no imposed boundary condition. The stiffness matrix \( \mat{A} \)
however will be singular, and have a zero-eigenvalue of multiplicity one. 

In addition to this, the computation of the smallest eigenvalue of a matrix is
a numerically unstable procedure. We will therefore \emph{estimate} the
condition numbers as follows:
\begin{align}
\cond{\mat{M}} \approx \frac{|\lambda_1(\mat{M})|}{|\lambda_n(\mat{M})|}\,,\text{ and } \,
\cond{\mat{A}} \approx \frac{|\lambda_1(\mat{A})|}{|\lambda_{n-1}(\mat{A})|},
\end{align}
using the second-smallest eigenvalue for the stiffness matrix.

\begin{figure}[htpb]
\centering
\subfloat[LRB]{\includegraphics[width=0.7\linewidth]{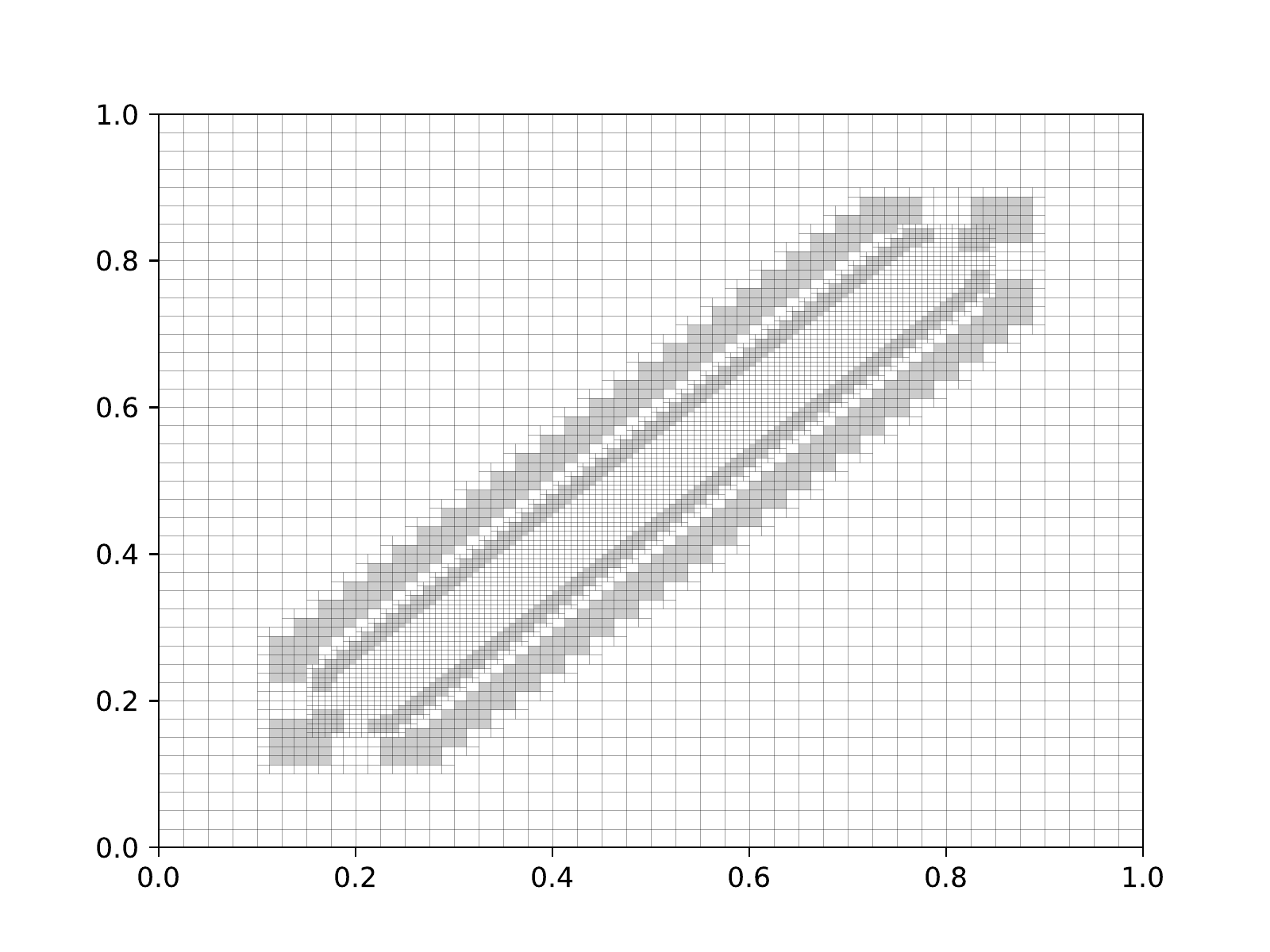}\label{sub:overloading_lrb_diag}}

\subfloat[THB]{\includegraphics[width=0.7\linewidth]{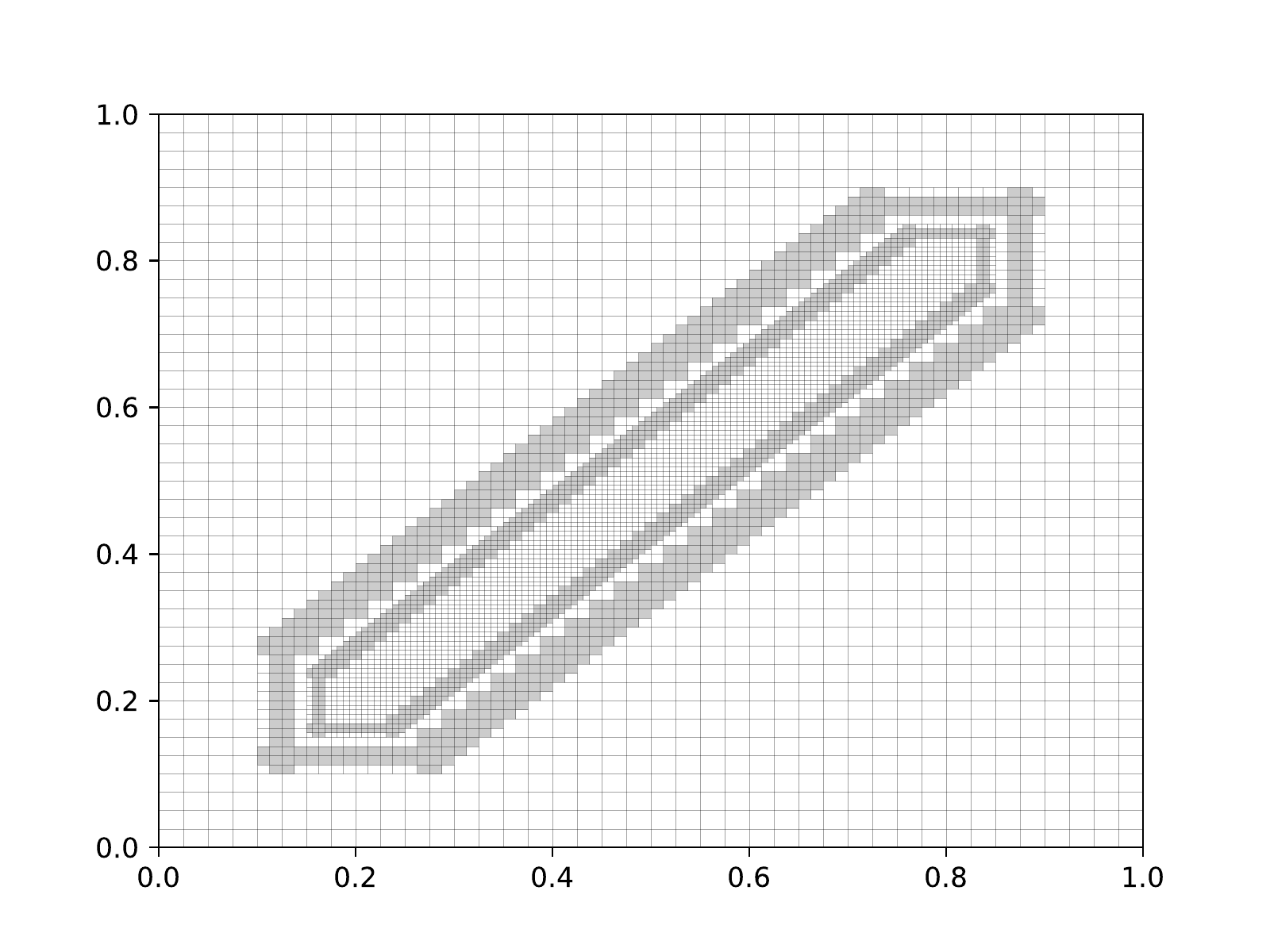}\label{sub:overloading_thb_diag}}

\subfloat[LRBNO]{\includegraphics[width=0.7\linewidth]{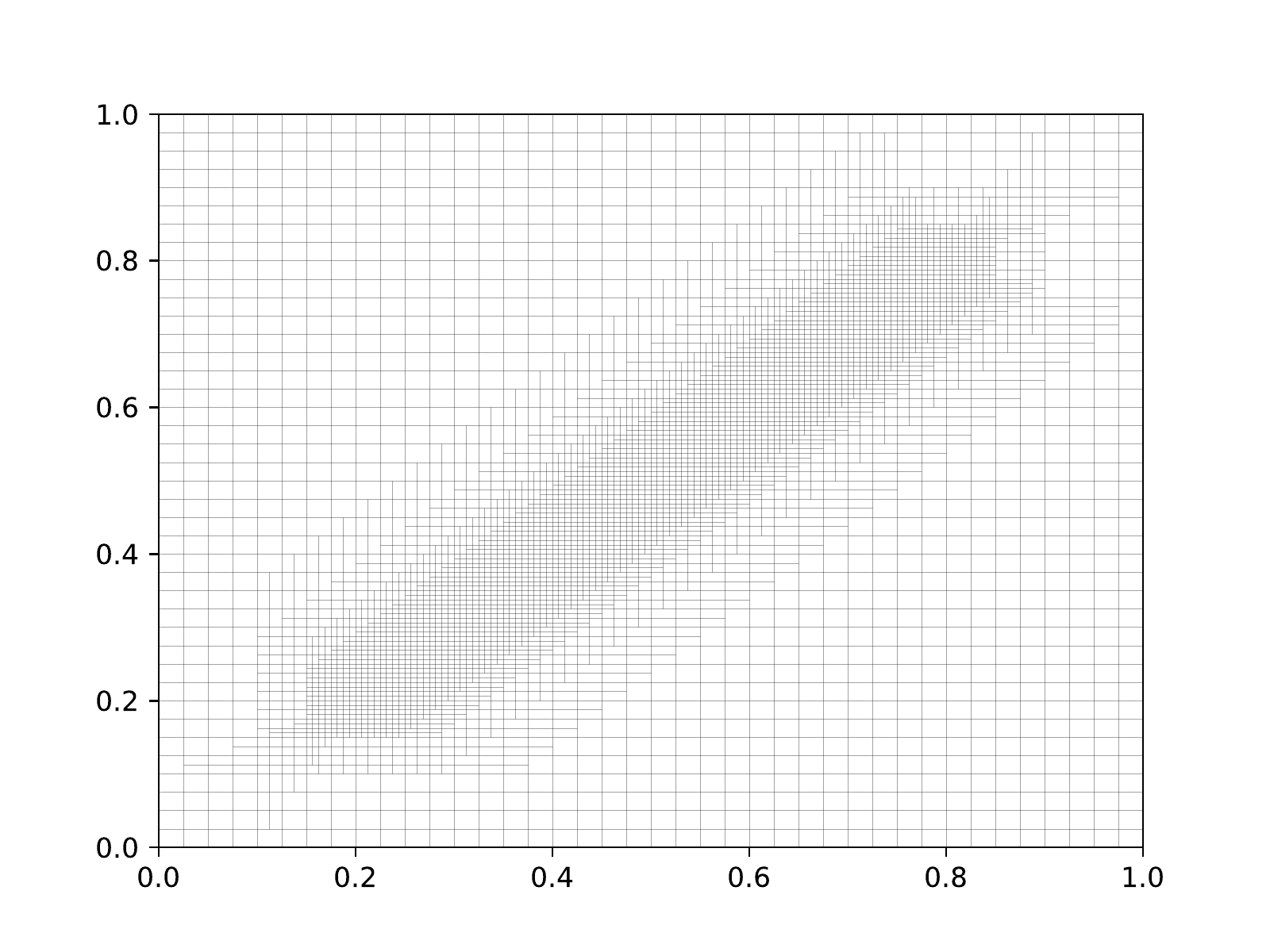}\label{sub:overloading_lrb_modified_diag}}

\caption[Overloading patterns for THB and unmodified LRB and modified LRB
(diagonal)]{The overloading patterns on a hierarchical mesh with three levels
of diagonal refinement. In this case, we see in greater effect the
behaviour of LRB over convex corners. Here the difference in overloading
between THB and LRB are smaller, as opposed to the central refinement case,
due to the high
number of corners relative to the length of the sides of the refined
levels. By using a one-directional meshline extension along the diagonal,
and extensions similar to the central-refinement case, we may completely
remove overloading.}
\label{fig:overloading_example_diag}
\end{figure}

\begin{figure}[htpb]
\centering

\includegraphics[width=1\linewidth]{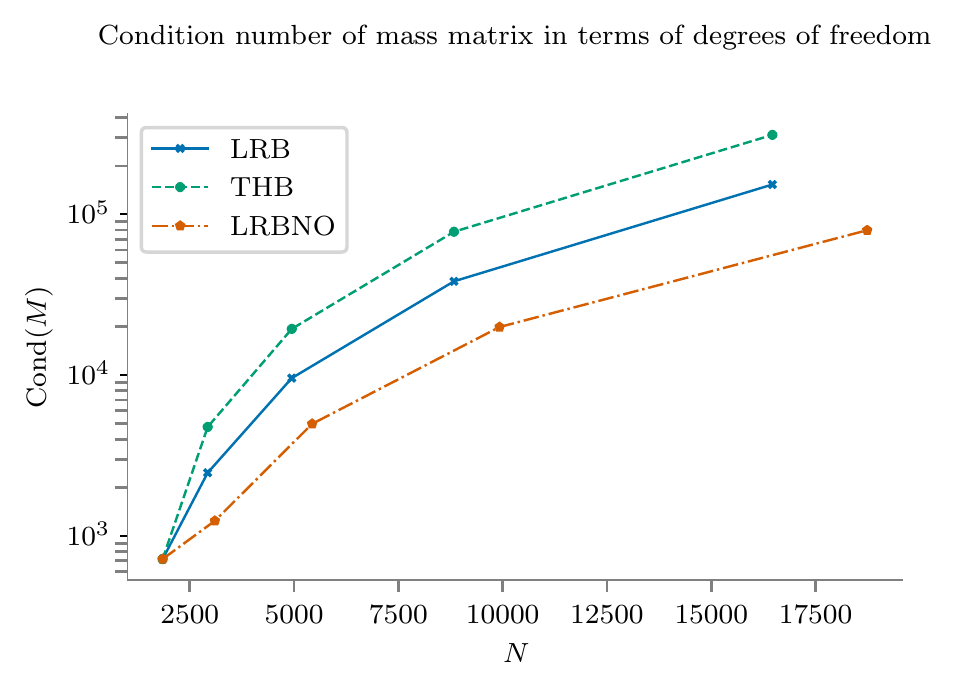}

\caption{Condition numbers for mass matrices over a diagonally refined
hierarchical mesh for four levels of refinement.  There is one data point for each method at each refinement level. The first point is the same for all methods as the all methods start from the same tensor product spline space. \( N \)
denotes the number of degrees of freedom in the corresponding spline space. The none overloaded LRBNO mesh has clearly the smallest condition numbers.}
\label{fig:diagonal_refinement_condition_mass}
\end{figure}

\begin{figure}[htpb]
\centering

\includegraphics[width=1\linewidth]{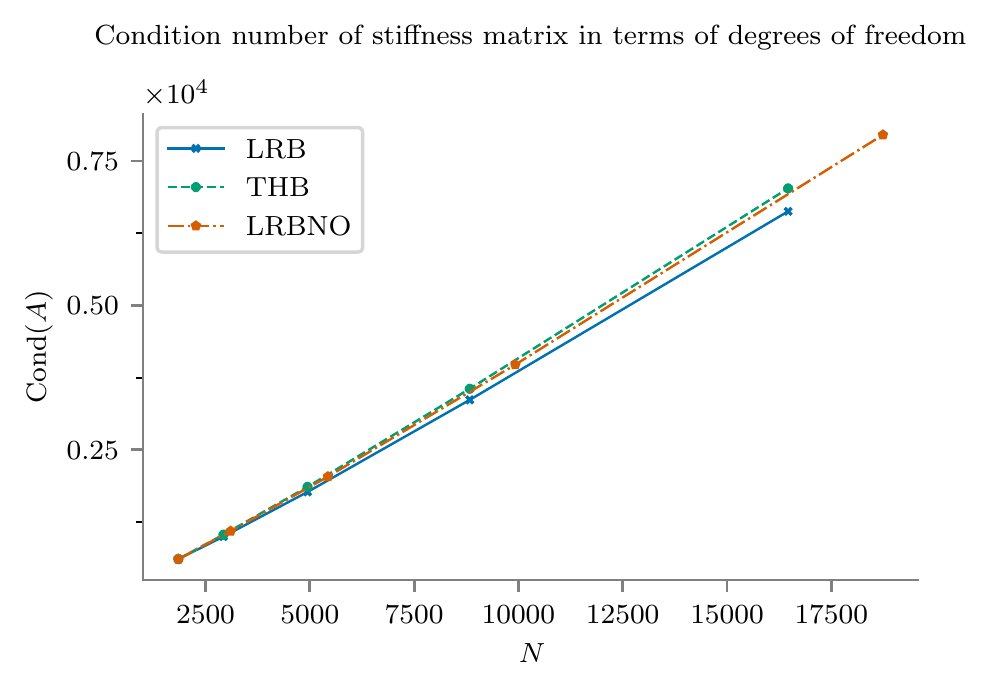}

\caption{The condition numbers for stiffness matrices over a
diagonally refined hierarchical mesh for four levels of refinement. In the
figures \( N \) denotes the number of degrees of freedom in the corresponding
space. All methods are similar in behaviour with respect to condition numbers as a function of degrees of freedom.}
\label{fig:diagonal_refinement_condition_stiffness}
\end{figure}

\subsection{Numerical results}

Below we present  results of the numerical simulations. As LRBNO generates higher dimensional spline spaces than THB and LRB we  plot  the condition numbers as a function of the
degrees of freedom. Just plotting the condition
numbers  as a function of the levels provides less information. By including the dimension of the spline space we obtain a clearer distinction between the methods.

\subsubsection{Central refinement}

We assemble the stiffness- and mass-matrices on a sequence of meshes
corresponding to central refinement, shown at the third refinement in
\cref{fig:overloading_example,fig:overloading_example_modified}. The results
are shown in
\cref{fig:central_refinement_condition_mass,fig:central_refinement_condition_stiffness}.
Start by noting that for the mass matrix, THB performs better than LRB with no
modifications, while for the stiffness matrices, the two methods are comparable
with LRB having a slight advantage. The number of degrees of freedom are the
same. By locally modifying the mesh, as is the case for LRBNO and T-LRBNO, we
see that the number of degrees of freedom goes up, as expected. The condition
number per degree of freedom is smallest for T-LRBNO.

\subsubsection{Diagonal refinements}

We now consider the case of diagonal refinement. Again, we use the same
hierarchical mesh for LRB and THB. We will only consider one sequence of meshes
with local modifications. In the diagonal refinement setting, the corners of
the refined region are sufficiently close to each other so that we need to make
a decision on which direction to refine in. The diagonally refined mesh is not
compatible with a T-spline type mesh, and will therefore not be taken into
consideration here.

We assembly stiffness- and
mass-matrices on the meshes displayed in \cref{fig:overloading_example_diag}.
The results are shown in \cref{fig:diagonal_refinement_condition_mass,fig:diagonal_refinement_condition_stiffness}. Note that
for the diagonal refinement, the number of degrees of freedom generated when
removing overloading, shown in the mesh in \cref{sub:overloading_lrb_modified_diag}, is a fair bit
larger than the unmodified counterparts. Despite this, LRBNO outperforms THB
and LRB by a significant amount when it comes to the mass matrix. The
conditioning of the stiffness matrix on the other hand grows approximately
linearly with the number of degrees of freedom, and no significant effect of
the overload-reduction can be seen.

%% file: conclusion.tex
\section{Conclusion}
\label{sec:conclusion}

We have addressed differences and similarities of Truncated Hierarchical
B-splines (THB) and Locally Refined B-splines (LRB) on similar hierarchical
meshes. The overall conclusion is that there are no big differences between
the methods with respect to condition numbers of mass and stiffness matrices
for the example meshes addressed. 
\begin{itemize}

\item When THB and LRB are run on identical meshes THB has better conditions
  numbers for the mass matrix except for most complex example run, the diagonal
  example in Figures \ref {fig:diagonal_refinement_condition_mass} and
  \ref{fig:overloading_example_diag}. The behaviour of the stiffness matrix is
  very similar for both methods.

\item When making a mesh for LRB that has no overloading the condition numbers
  for the mass matrix of LRB are smaller than those of THB, with condition
  numbers of stiffness being similar. It should be noted that using meshes for
  LRB that has no overloading guarantees that the B-splines generated are
  linearly independent, and that the number of B-splines covering an element is
  the minimal  needed for spanning the polynomial space over the element. For
  hierarchical meshes of bi-degree less than $(4,4)$ there is always linear
  independence in the set of LR B-splines generated,. For bi-degree $(4,4)$ and
  higher linear dependence can occur in very special configurations when the
  elements outside two opposing concave corners of a refinement region is
  covered by the same B-spline from a cruder level. This happens for bi-degree
  (4,4) when a refinement region is split if just one element from the cruder
  level is not refined, e.g., the refinement region is locally very narrow.
\end{itemize}

When trying to represent hierarchical refinements using T-splines as in Figure
\ref{fig:t_mesh} there is a region of one directional refinement of length two
just outside the boundary of the refinement region. This gives a smoother
transition between refinement levels that can also be replicated by LRB. The
results in Figures \ref{fig:central_refinement_condition_mass} and
\ref{fig:central_refinement_condition_stiffness} show a better behaviour than
going directly from one refinement level to the next.  Having such an
intermediate level of refinement if possible is advantageous.  However, in
situations such as the diagonal refinement in Figure
\ref{fig:overloading_example_diag} this is not possible.

Most often THB is described as based on dyadic sequences of grids determined by
scaled lattices over which uniform B-spline spaces are defined. This implies
that there is single knot multiplicity along domain boundaries. However,
variants of THB are published \cite{GIANNELLI2016337} where open knots are used
along the domain boundary. In Section \ref{sec:boundary_multiplicity} we have
shown that open knot vectors are preferable, not only with respect to
simplified interpolation of boundary conditions, but also to avoid that the
condition number of the mass matrix is biased by the boundary B-splines. As we
see the same effect for LR B-splines we have a strong recommendation that open
knot vectors are used for locally refined splines, rather than single knot
multiplicity at domain boundaries.